\documentclass[12 pt]{article}
\usepackage{amssymb,amsmath,amsfonts,amsthm}
\renewcommand{\P}{\ensuremath{\mathbb {P}}}

\newtheorem {Theorem}{Theorem}[section]

\newtheorem {Corollary}[Theorem]{Corollary}

\theoremstyle{definition}

\newtheorem{Remark}{Remark}[section]

\newcommand\beq{\begin{equation}}
\newcommand\eeq{\end{equation}}

\begin{document}
\title{Behavior of the empirical Wasserstein distance in ${\mathbb R}^d$ under moment conditions}

\date{}

\author{ J\'er\^ome Dedecker$^a$ and Florence Merlev\`ede$^b$}

\maketitle

{\abstract{We establish some deviation inequalities, moment bounds and almost sure results for the Wasserstein distance of order $p\in [1, \infty)$ between 
the empirical measure of independent and identically distributed ${\mathbb R}^d$-valued random variables and the common  distribution 
of the variables. We only assume the existence of a (strong or weak) moment  of order $rp$ for some $r>1$, and we discuss the optimality of the bounds.} }

\medskip


\medskip


\bigskip

\small{

\noindent $^a$ Universit\'e Paris Descartes, Sorbonne Paris Cit\'e,  Laboratoire MAP5 (UMR 8145).\\
Email: jerome.dedecker@parisdescartes.fr \bigskip

\noindent $^b$ Universit\'{e} Paris-Est, LAMA (UMR 8050), UPEM, CNRS, UPEC. \\
Email: florence.merlevede@u-pem.fr}

\section{Introduction and notations}

\setcounter{equation}{0}
We begin with some notations, that will be used all along the paper. 
Let $X_1, \ldots, X_n$ be $n$ independent and identically distributed (i.i.d.) random variables with values in ${\mathbb R}^d$ ($d \geq 1$), with
 common distribution $\mu$. Let $\mu_n$ be the empirical distribution 
 of the $X_i$'s, that is 
 $$
   \mu_n= \frac 1 n \sum_{k=1}^n  \delta_{X_k} \, .
 $$
 
 Let  $X$ denote a random variable with distribution $\mu$. 
For any $x \in {\mathbb R}^d$, let $|x|= \max \{ |x_1|, \ldots, |x_d| \}$. Define then the tail  of the distribution $\mu$ by
$$
  H(t)= {\mathbb P}(|X|>t) = \mu( \{ \text{$x  \in {\mathbb R}^d$ such that  $|x|>t$}\}) \, .
$$
As usual, for any $q \geq 1$, the weak moment of order $q$ of the random variable $X$ is defined by  $$\|X\|_{q, w}^q:= \sup_{t >0} t^q H(t) \, ,$$ 
and the strong moment of order $q \geq 1$ is defined by 
$$
 \|X\|_{q}^q = {\mathbb E} \left( |X|^q\right )=  q\int_0^\infty t^{q-1} H(t) dt  \, .
$$

For $p \geq 1$, the Wasserstein distance between two probability measures $\nu_1, \nu_2$  on $({\mathbb R}^d, {\mathcal B}({\mathbb R}^d))$ 
is defined by 
$$
   W_p^p(\nu_1, \nu_2) = \inf_{ \pi \in \Pi(\nu_1, \nu_2) } \int_{{\mathbb R}^d \times {\mathbb R}^d}  |x-y|_2^p  \ \pi(dx,dy) \, ,
$$
where $|\cdot |_2$ is the euclidean norm on ${\mathbb R}^d$ and $\Pi(\nu_1, \nu_2)$ is the set of probability measures on the product space 
$({\mathbb R}^d \times {\mathbb R}^d , {\mathcal B}({\mathbb R}^d) \otimes   {\mathcal B}({\mathbb R}^d))$  with margins $\nu_1$ and $\nu_2$.

In this paper, we prove deviation inequalities, moment inequalities and almost sure results  for the quantity $W_p(\mu_n, \mu)$, 
when $X$ has a weak or strong moment of order $rp$ for $r>1$. 
As in \cite{FG}, the upper bounds will be different according as $p>d\min \{(r-1)/r, 1/2\}$ (small dimension case) or $p<d \min \{(r-1)/r, 1/2\}$ (large dimension case). 
Most of the proofs are based on Lemma 6 in \cite{FG} (see the inequality \eqref{depart} in Section \ref{Proofs}), which may be seen as an extension of \`Ebralidze's inequality \cite{E} to the case $d>1$. 
Hence we shall use the same approach as in \cite{DM}, where we combined \`Ebralidze's inequality with truncation arguments to get moment bounds
for $W_p(\mu_n, \mu)$ when $d=1$. 

There are many ways to see that the upper bounds obtained in the present paper are optimal in some sense,
by considering the special cases $d=1$, $p=1$, $p=2$,  or by following the general discussion in \cite{FG}, and we shall make some comments about this question all along the paper.
However, the optimality for large $d$ is only a kind of minimax optimality: one can see that the rates are exact for compactly supported measures
which are not singular with respect to the Lebesgue measure on ${\mathbb R}^d$ (by using, for instance, Theorem 2 in \cite{DSS}).  

In fact, since the rates depend on the dimension $d$, it is easy to see 
that they cannot be optimal for all measures: for instance the rates will be faster as announced if the measure $\mu$ is supported on a linear subspace of ${\mathbb R}^d$ with dimension strictly
less than $d$. This is of course not the end of the story, and the problem can be formulated in the general context of metric spaces $(X, \delta)$. For instance, for compactly supported measures, Boissard and Le Gouic \cite{BLG} proved that the rates of convergence depend on the behavior of the  metric entropy of the support of 
$\mu$  (with an extension to non-compact support in their Corollary 1.3). In the same context, Bach and Weed \cite{BW} obtain sharper results by generalizing some ideas
going back to Dudley (\cite{D}, case $p=1$). They introduce the notion of Wasserstein dimension $d^*_p(\mu)$  of the measure $\mu$, and prove that 
$n^{p/s}{\mathbb E}(W^p_p(\mu_n, \mu))$ converges to 0  for any $s>d^*_p(\mu)$ (with sharp lower bounds in most cases). 

Note that our context and that of Bach and Weed are clearly distinct:   we consider measures on ${\mathbb R}^d$ having only a finite moment of order $rp$ for $r>1$, 
while they consider measures on compact metric spaces. However, the Wasserstein dimension is well defined for any probability  measure (thanks to Prohorov's theorem), and 
some arguments in \cite{BW} are common with \cite{DSS} and \cite{FG}. A reasonable question is then:  in the case of a singular measures on  ${\mathbb R}^d$, are the results of the present paper still valid if we replace the dimension
 $d$ by any $d' \in (d^*_p(\mu), d]$?
 
 \smallskip
 
 The paper is organized as follows: in Section \ref{di} we state some deviations inequalities for $W_p(\mu_n, \mu)$ under weak moment assumptions.
 In Section \ref{LMD} we bound up the probability of large and moderate deviations. In Section \ref{ASR} we present some almost sure results,
 and in Section \ref{MI} we give some upper bounds for the moments of order $r$ of $W_p(\mu_n, \mu)$ (von Bahr-Esseen and Rosenthal type bounds)
 under strong moment assumptions. The proofs are given in Section \ref{Proofs}.
 
 \smallskip
 
 All along the paper, we shall use the notation $f(n, \mu, x) \ll g(n, \mu, x)$, which means that there exists a positive constant $C$, not depending 
 on $n, \mu, x$ such that $f(n, \mu, x) \leq C g(n, \mu, x)$ for all positive integer $n$ and all positive real $x$.

\section{Deviation inequalities under weak moments conditions}\label{di}

\setcounter{equation}{0}

In this section, we give some upper bound for the quantity 
 ${\mathbb P}( W_p^p(\mu_n, \mu) >x) $ when the random variables $X_i$ have a weak moment of order $rp$ for some $r>1$. 
 We first consider the case where $r \in (1,2)$. 

\begin{Theorem}\label{weakVBE}
If $\|X\|_{rp, w} < \infty$ for some $r \in (1,2)$, then 

\[
 {\mathbb P}( W_p^p(\mu_n, \mu) >x)  \ll  \left\{ \begin{aligned}
 & \frac{\|X\|_{rp, w}^{rp}}{x^r n^{r-1}}  & \text{ if } & p > d(r-1)/r\\
 &   \frac{\|X\|_{rp, w}^{rp}(\log n)^r}{x^r n^{r-1}}  \left ( 1 + \log_{+} \left( \frac{x^{1/p} n^{r/(dr-d))}}{\|X\|_{rp, w}} \right) \right)^r  & \text{ if } & p  = d(r-1)/r\\
 &   \frac{\|X\|_{rp, w}^{rp}}{x^r n^{rp/d}}    & \text{ if } & p  \in [1, d(r-1)/r)
    \end{aligned}
    \right.
\]
for any $x>0$, where $\log_{+}(x)= \max \{0, \log x  \}$.
\end{Theorem}

\begin{Remark}
As will be clear from the proof, the upper bounds of Theorem \ref{weakVBE} still hold if the quantity
${\mathbb P}( W_p^p(\mu_n, \mu) >x)$ is replaced by its maximal version
$$
{\mathbb P}\left ( \max_{1 \leq k \leq n}kW_p^p(\mu_k, \mu) >nx \right ) \, . 
$$
Since $\| W_p^p(\mu_n, \mu)  \|_1 \leq (r /(r-1)) \,  \|W_p^p(\mu_n, \mu)\|_{r, w}$,  according to the discussion after Theorem 1 in \cite{FG}, if $p \neq d(r-1)/r$, one can always find some measure $\mu$ 
for which the rates of Theorem \ref{weakVBE} are reached (see example (e)  in \cite{FG} for $ p > d(r-1)/r$ and example (c) 
 in \cite{FG} for $ p < d(r-1)/r$).

\end{Remark}

\smallskip

We now consider the case where $r>2$. We follow the approach of Fournier and Guillin \cite{FG}, but we use a different 
upper bound for the quantity controlled in their Lemma 13 (see the proof of Theorem \ref{F-Nag} for more details).

\begin{Theorem}\label{F-Nag}
If $\|X\|_{rp, w} < \infty$ for some $r \in (2, \infty)$, then  for 
 any $q>r$,
\[
 {\mathbb P}( W_p^p(\mu_n, \mu) >x)  \ll  
 a \left (n, \frac{x}{ \|X\|_{rp, w}^{p}} \right ) + \frac{\|X\|_{rp, w}^{rp}}{x^r n^{r-1}}  + \frac{1}{x^q n^{q/2}} 
\left ( \int_0^\infty t^{p-1} \sqrt {H(t)} dt \right )^q \, ,
\]
for any $x>0$,
where
\[
 a(n,x) = C  \left\{ \begin{aligned}
 & \exp (-c n x^2){\mathbf 1}_{x \leq A} & \text{ if } & p > d/2 \\
 &\exp (-c n  ( x / \log ( 2 + x^{-1} ) )^2 )
 {\mathbf 1}_{x \leq A} & \text{ if } & p = d/2 \\
 &   \exp (-c n x^{d/p} ){\mathbf 1}_{x \leq A} & \text{ if } & p  \in [1, d/2) 
    \end{aligned}
    \right.
\]
for some positive constants $C,c$ depending only on $p,d$, and a positive  constant $A$ depending only on  $p,d,r$. 
\end{Theorem}

\begin{Remark}
Let us compare our inequality with that of Theorem 2  of Fournier and Guillin \cite{FG} (under the moment condition (3) in \cite{FG}). 
We first note that the inequality in \cite{FG} is stated under a strong moment of order $rp$ for $r>2$, but their proof
works also under a weak moment of order $rp$. Hence, under the assumptions of our Theorem \ref{F-Nag}, Fournier and Guillin obtained the bound
(we assume here  that  $\|X\|_{rp, w}=1$ for the sake of simplicity):
\begin{equation}\label{FGberk}
{\mathbb P}( W_p^p(\mu_n, \mu) >x)  \ll  
 a \left (n, x \right ) + \frac{n}{(nx)^{(rp- \varepsilon)/p}} \, ,
\end{equation}
for any $\varepsilon>0$ (the constant implicitly  involved in the inequality depending on $\varepsilon$). In particular, one cannot infer from \eqref{FGberk}
that 
$$
\limsup_{n \rightarrow \infty} 
n^{ r -1} {\mathbb P}\left ( W_p^p(\mu_n, \mu) > x
\right )  \ll  \frac{\|X\|_{rp, w}^{rp}}{x^r} \, ,
$$
which follows from our Theorem \ref{F-Nag}. 
\end{Remark}

\section{Large and moderate deviations}\label{LMD}
\setcounter{equation}{0}

We consider here the probability of moderate deviations, that 
is 
$$
 {\mathbb P}\left ( W_p^p(\mu_n, \mu) >\frac{x}{n^{1- \alpha}}
\right ) \, ,
$$
for $\alpha \leq 1$ in a certain range and $x>0$. As usual, the case $\alpha =1$ is the probability 
of large deviations. 

As for partial sums, we shall establish two type of results, under weak moment conditions  or under strong
moment conditions. If the random variables have a weak moment of order $rp$ for some $r>1$, the results  of Subsection \ref{wm} are immediate corollaries of the  theorems of the preceding section. 
On the contrary, the Baum-Katz type results of Subsection \ref{bk} cannot be derived from the results of Section \ref{di} and will be proved in Subsection \ref{proofbk}.

\subsection{Weak moments}\label{wm}
As a consequence of Theorem \ref{weakVBE}, we obtain the following corollary.

\begin{Corollary}\label{moderate1}
 If $\|X\|_{rp, w} < \infty$ for some $r \in (1,2)$, then,
\begin{itemize}
\item If $p>d(r-1)/r$ and $1/r \leq \alpha \leq 1$, 
$$
\limsup_{n \rightarrow \infty} 
n^{\alpha r -1} {\mathbb P}\left ( W_p^p(\mu_n, \mu) >\frac{x}{n^{1- \alpha}}
\right )  \ll  \frac{\|X\|_{rp, w}^{rp}}{x^r} \, .
$$
\item  If $p=d(r-1)/r$ and $1/r <  \alpha \leq 1$, 
$$
\limsup_{n \rightarrow \infty} 
\frac{n^{\alpha r -1}}{(\log n)^{2r}} {\mathbb P}\left ( W_p^p(\mu_n, \mu) >\frac{x}{n^{1- \alpha}}
\right )  \ll  \frac{\|X\|_{rp, w}^{rp}}{x^r} \, .
$$
\item If $p \in [1, d(r-1)/r)$ and 
$
 (d-p)/d \leq \alpha \leq 1
$,
$$
\limsup_{n \rightarrow \infty} 
 n^{(pr-(1-\alpha)rd)/d}\,  {\mathbb P}\left ( W_p^p(\mu_n, \mu) >\frac{x}{n^{1- \alpha}}
\right )  \ll  \frac{\|X\|_{rp, w}^{rp}}{x^r} \, .
$$
\end{itemize}
\end{Corollary}

\begin{Remark}
Let us comment on the case $p=1, d=1$. In that case, del Barrio et al. \cite{dBGM} proved that, for $\beta \in (1,2)$,
$ n^{(\beta-1)/\beta} W_1( \mu_n, \mu)$ is stochastically bounded if and only if 
$\|X\|_{\beta, w} < \infty$ (see their Theorem 2.2). This is consistent with the first inequality of Corollary  \ref{moderate1} applied with $r=\beta$ and $\alpha=1/r$. 
\end{Remark}

\begin{Remark}
Let us now comment on the case $p=2, d=1$. In that case del Barrio et al. \cite{dBGU} proved that, if the distribution function $F$ of $X$ is twice differentiable
and if $F' \circ F^{-1}$ is a regularly varying function in the neighborhood of 0 and 1, then there exists a sequence of positive numbers $v_n$ tending to $\infty$ as $n \rightarrow \infty$, such that 
$v_n W_2^2(\mu_n, \mu)$ converges in distribution to a non degenerate distribution. For instance, it follows from their Theorem 4.7 that, if $X$ is a positive random variable, $F$ is twice differentiable and $F(t)= (1-t^{-\beta})$ for any $t >t_0$ and some $\beta >2$, then $n^{(\beta-2)/\beta}W_2^2(\mu_n, \mu)$ converges in distribution to a non degenerate distribution.
In that case, there is a weak moment of order $\beta$, and,  for $\beta \in (2,4)$, the first inequality of Corollary \ref{moderate1} applied with $r=\beta/2$ and $\alpha=1/r$ gives
$$
\limsup_{n \rightarrow \infty} 
 {\mathbb P}\left ( n^{(\beta-2)/\beta} W_2^2(\mu_n, \mu) >x
\right )  \ll  \frac{\|X\|_{\beta, w}^{\beta}}{x^{\beta/2}} \, .
$$
Hence, in the case where $\beta \in (2,4)$, our result is consistent with that given in \cite{dBGU}, and holds without assuming any regularity on $F$. 
\end{Remark}

As a consequence of Theorem \ref{F-Nag}, we obtain the following corollary.

\begin{Corollary}
 If $\|X\|_{rp, w} < \infty$ for some $r \in (2, \infty)$, then, for  any 
$$
 \alpha \in \left ( \max \left ( \frac 1 2, \frac {d-p}{d} \right ) ,  1
 \right ] \, ,
$$
$$
\limsup_{n \rightarrow \infty} 
n^{\alpha r -1} {\mathbb P}\left ( W_p^p(\mu_n, \mu) >\frac{x}{n^{1- \alpha}}
\right )  \ll  \frac{\|X\|_{rp, w}^{rp}}{x^r} \, .
$$
\end{Corollary}

\subsection{Baum-Katz type results} \label{bk}

In this subsection, we shall prove some deviation results in the spirit of Baum and Katz \cite{BK}. Recall that, for partial sums
$S_n=Y_1+\cdots + Y_n$ of i.i.d real-valued random variables such that $\|Y_1\|_r <\infty$  for some $r>1$ and ${\mathbb E}(Y_1)=0$, one has:
for any  $\alpha >1/2$ such that $1/r \leq \alpha \leq 1$,   and any $x>0$, 
$$
\sum_{n=1}^\infty n^{\alpha r-2} {\mathbb P} \left ( \max_{1 \leq k \leq n} |S_k| > n^{\alpha} x
\right )  < \infty  \, .
$$

We first consider the case where the variables have a strong moment of order $rp$ for $r \in (1,2)$. 

\begin{Theorem}\label{BK1} If $\|X\|_{rp} < \infty$ for some $r \in (1,2)$, then, for any $x>0$,
\begin{itemize}
\item If $p>d(r-1)/r$ and $1/r \leq \alpha \leq 1$, 
$$
\sum_{n=1}^{\infty} 
n^{\alpha r -2} {\mathbb P}\left ( \max_{1 \leq k \leq n} k W_p^p(\mu_k, \mu) > n^{\alpha} x
\right )  < \infty  \, .
$$
\item If $p \in [1, d(r-1)/r)$ and 
$
 \alpha \in \left ( (d-p)/d ,  1
 \right ] \, ,
$
$$
\sum_{n=1}^{\infty} 
n^{(pr-(1-\alpha)rd -d)/d} \, 
{\mathbb P}\left ( \max_{1 \leq k \leq n} k W_p^p(\mu_k, \mu) > n^{\alpha} x
\right )  < \infty   \, .
$$
\item If $p \in [1, d(r-1)/r)$, 
$$
\sum_{n=1}^{\infty} 
\frac 1 n  \, 
{\mathbb P}\left ( \max_{1 \leq k \leq n} k W_p^p(\mu_k, \mu) > n^{(d-p)/d} (\log n)^{1/r} x
\right )  < \infty   \, .
$$
\end{itemize}
\end{Theorem}

\begin{Remark}
Our proof does not allow to deal with the case where $p=d(r-1)/r$. As an interesting consequence of  Theorem \ref{BK1}, we shall obtain almost sure convergence rates 
for the sequence $W_p^p(\mu_n, \mu)$ (see Corollary \ref{asrates} of the next section). 
\end{Remark}

We now consider the case where the variables have a strong moment of order $rp$ for $r  >2$. 

\begin{Theorem}\label{BK2}
If $\|X\|_{rp} < \infty$ for some $r \in (2, \infty)$, then, for any $x>0$ and any
$$
 \alpha \in \left ( \max \left ( \frac 1 2, \frac {d-p}{d} \right ) ,  1
 \right ] \, ,
$$
$$
\sum_{n=1}^{\infty} 
n^{\alpha r -2} {\mathbb P}\left ( \max_{1 \leq k \leq n} k W_p^p(\mu_k, \mu) > n^{\alpha} x
\right )  < \infty  \, .
$$
\end{Theorem}

\section{Almost sure results}\label{ASR}
\setcounter{equation}{0}

Using well known arguments, we derive from Theorem \ref{BK1} the following almost sure rates of convergence for the sequence 
$W_p^p(\mu_n, \mu)$
(taking  $\alpha= 1/r$ in the case where $p> d(r-1)/r$, and applying the third item in the case where $p< d(r-1)/r$). 

\begin{Corollary}\label{asrates}
If $\|X\|_{rp} < \infty$ for some $r \in (1,2)$, then
\begin{itemize}
\item If $p>d(r-1)/r$,
$$
\lim_{n \rightarrow \infty} n^{(r-1)/r} W_p^p ( \mu_n, \mu)=0 \, \text{ a.s. }
$$
\item If $p \in [1, d(r-1)/r)$, 
$$
\lim_{n \rightarrow \infty} \frac{ n^{p/d}} {(\log n)^{1/r}}W_p^p ( \mu_n, \mu)=0 \, \text{ a.s. }
$$
\end{itemize}
\end{Corollary}

\begin{Remark}
Let us comment on these almost sure results in the case where $p=1$ and $d < r/(r-1) $. Recall the dual expression of $W_1(\mu_n, \mu)$:
\begin{equation}\label{dual}
W_1(\mu_n, \mu)=  \sup_{f \in \Lambda_1} \left | \frac 1 n  \sum_{k=1}^n \left  ( f(X_k) - \mu(f)  \right ) \right | 
\end{equation}
where $\Lambda_1$ is the the set of functions $f$ such that $|f(x)-f(y)| \leq |x-y|_2$.
Since the function $g: x \mapsto |x|_2$ belongs to $\Lambda_1$, we get
\begin{equation*}\label{clear}
   W_1(\mu_n, \mu) \geq \frac 1 n \left | \sum_{k=1}^n \left ( \, |X_k|_2- {\mathbb E} ( |X_k|_2) \, \right) \right | \, . 
\end{equation*}
Now, by the classical   Marcinkiewicz-Zygmund theorem (see \cite{MZ}) for i.i.d. random variables, we know that 
$$
\lim_{n \rightarrow \infty}   \frac {n^{(r-1)/r}} {n} \left | \sum_{k=1}^n \left ( \, |X_k|_2- {\mathbb E}  (|X_k|_2) \, \right) \right |  = 0 \, \text{ a.s. }
$$
if and only if $\|X\|_{r} < \infty$. It follows that, for $p=1$,  the rates given in Corollary \ref{asrates} are optimal in the case where $d<r/(r-1)$. 
\end{Remark}

\smallskip

We now give some almost sure rates of convergence in the case where $ \int_0^{\infty} t^{p-1} \sqrt{ H(t) } dt < \infty$. 
Note that this condition is a bit more restrictive than $\|X\|_{2p} < \infty$ (but is satisfied, for instance, 
if  $\|X\|_{rp} < \infty$ for some $r>2$). 

\begin{Theorem}\label{LILBornee} 
Assume that  $ \int_0^{\infty} t^{p-1} \sqrt{ H(t) } dt < \infty$. 
\begin{itemize}
\item If $p>d/2$,  there exists an universal positive constant $C$ depending only on  $(p,d)$ such that 
\[
 \limsup_{n \rightarrow \infty}  \sqrt{\frac{ n }{\log \log n} }  W_p^p(\mu_n, \mu)  \leq  C   \int_0^{\infty} t^{p-1} \sqrt{ H(t) } dt \quad a.s.
\]
\item If $p\in [1,d/2)$,  there exists an universal positive constant $C$ depending only on  $(p,d)$ such that 
\[
 \limsup_{n \rightarrow \infty}   \left ( \frac{ n }{\log \log n} \right )^{p/d}   W_p^p(\mu_n, \mu)  \leq  C   \int_0^{\infty} t^{p-1} \sqrt{ H(t) } dt \quad a.s.
\]
\end{itemize}
\end{Theorem}

\begin{Remark} In the case $p>d/2$, 
the rate $\sqrt{n /\log \log n}$ has been obtained recently by Dolera and Reggazini (\cite{DR}, Theorem 2.3)  under the more restrictive condition $\|X\|_{rp} < \infty$ for some $r>2$.
\end{Remark}

\begin{Remark} 
In the case $p=1,d=1$, it follows from the central limit theorem for $W_1(\mu_n, \mu)$ (see \cite{dBGM}) and from Theorem 10.12 in \cite{LT} that
the sequence $( \sqrt{n /\log \log n} \ W_1(\mu_n, \mu))_{n \geq 0}$ is almost surely relatively compact if $\int_0^{\infty}  \sqrt{ H(t) } dt < \infty$, 
which is consistent with the first item of Theorem \ref{LILBornee}. 
\end{Remark}

\begin{Remark}
For $p=1$, concerning the rate of Corollary \ref{asrates} when $d > r/(r-1) $ or the rate of Theorem \ref{LILBornee} when $d> 2$, the situation is not as clear as in the small dimension case.  
According to Talagrand \cite{T}, if $d >2 $ and $\mu$ is the uniform measure on $[0,1]^d$, $W_1(\mu_n, \mu)$  is, almost surely,   exactly of
order $n^{-1/d}$. More generally, let us recall a  result by Dobri\'c and Yukich \cite{DY}: if   $d >2$ and $\mu$   is compactly supported, then, almost surely,
\begin{equation}\label{DY}
c'(d) \int (f_\mu(x))^{(d-1)/d} \leq   \liminf_{n \rightarrow \infty}  n^{1/d} W_1(\mu_n, \mu) \leq  \limsup_{n \rightarrow \infty}  n^{1/d} W_1(\mu_n, \mu) \leq  c(d) \int (f_\mu(x))^{(d-1)/d}  
\end{equation}
where $c(d), c'(d)$ depend only on $d$, and $f_\mu$ is the density of the absolutely continuous part of $\mu$ (hence the limit is zero if $\mu$ is singular with respect to the Lebesgue measure
on ${\mathbb R}^d$). Actually, it was announced in \cite{DY} that $c'(d)= c(d)$, but a gap in the proof  has been pointed out in \cite{BB}. 
\end{Remark}

\begin{Remark}
If $p<d/2$, Barthe and Bordenave \cite{BB} (see their Theorem 2) proved that, almost surely, 
\begin{equation}\label{BB} 
\beta'_p(d) \int (f_\mu(x))^{(d-p)/d}  \leq 
 \liminf_{n \rightarrow \infty}  n^{p/d} W_p^p(\mu_n, \mu)  \leq 
 \limsup_{n \rightarrow \infty}  n^{p/d} W_p^p(\mu_n, \mu) \leq  \beta_p(d) \int (f_\mu(x))^{(d-p)/d}  
\end{equation}
provided $\|X\|_{rp} < \infty$ for some $r>4d/(d-2p)$, which is a generalization of \eqref{DY}. For  $p<d/2$, Theorem \ref{LILBornee}  is difficult to compare with   \eqref{BB}, because the results  do not hold under the same assumptions on $d$ and $H$. 
A  reasonable questions is: 
does \eqref{BB} hold if  $\int_0^{\infty}  t^{p-1} \sqrt{ H(t) } dt < \infty$  and $p<d/2$? 
\end{Remark}

\section{Moment inequalities}\label{MI}

\setcounter{equation}{0}

In this section, we give some upper bounds for the moments  $\left \| W_p^p(\mu_n, \mu)  \right \|_r^r$ when the variables have a strong
moment of order $rp$. 

As will be clear from the proofs, the maximal versions of these inequalities hold, namely: the quantity 
$\left \| W_p^p(\mu_n, \mu)  \right \|_r$ can be replaced by  $$ \frac {1}{n^r} \left \|  \max_{1 \leq k \leq n}  k W_p^p(\mu_k, \mu)  \right \|_r^r$$ in all the statements of this section.

\subsection{Moment of order $1$ and $2$}

\begin{Theorem}\label{order1}
Let $q \in (1,2]$. If $\|X\|_p < \infty$, then, for any $M>0$, 
\[
 \left \| W_p^p(\mu_n, \mu)  \right \|_1 \ll  \left\{ \begin{aligned}
 & 
 \int_0^{\infty} t^{p-1} H(t) {\bf 1}_{t>M}  \, dt 
 +
 \frac{1}{ n^{(q-1)/q} }
 \int_0^{\infty} t^{p-1} (H(t))^{1/q} {\bf 1}_{t \leq M} \,  dt 
 & \text{ if } & p > d(q-1)/q \\
 & \int_0^{\infty} t^{p-1} H(t) {\bf 1}_{t>M}  \, dt 
 +
 \frac{\log n}{n^{p/d}} 
 \int_0^{\infty} t^{p-1} (H(t))^{(d-p)/d} {\bf 1}_{t \leq M } \,  dt & \text{ if } & p = d(q-1)/q \\
 &   \int_0^{\infty} t^{p-1} H(t) {\bf 1}_{t>M}  \, dt 
 +
 \frac{1}{ n^{p/d}} 
 \int_0^{\infty} t^{p-1} (H(t))^{(d-p)/d} {\bf 1}_{t \leq M} \,  dt 
 & \text{ if } & p  \in [1, d(q-1)/q)
    \end{aligned}
    \right.
\]
where the constant implicitly involved does not depend on $M$.
\end{Theorem}
\begin{Remark}
In particular, if $H(t) \leq C t^{-p} (\log(1+t))^{-a}$ for some  $C>0, a>1$, then 
$$
\left  \| W_p^p(\mu_n, \mu)  \right \|_1=O\left(\frac{1}{(\log n)^{a-1}}\right )\, .
$$
\end{Remark}

\begin{Remark}
If $\|X\|_{rp, w} < \infty$ for $r \in (1,2)$ and $p \neq d(r-1)/r$, we easily infer  from Theorem \ref{order1} that
\[
 \left \| W_p^p(\mu_n, \mu)  \right \|_1  \ll  \left\{ \begin{aligned}
 & \frac{\|X\|_{rp, w}^{p}}{n^{(r-1)/r}}  & \text{ if } & p > d(r-1)/r \\
 &   \frac{\|X\|_{rp,w}^{p}}{n^{p/d}}   & \text{ if } & p  \in [1, d(r-1)/r) 
    \end{aligned}
    \right.
\]
which can also be deduced from Theorem \ref{weakVBE}. If $p=d(r-1)/r$, we get 
$$
 \left \| W_p^p(\mu_n, \mu)  \right \|_1 \ll \frac{\|X\|_{rp,w}^{p} (\log n)^{2}}{n^{p/d}}  \, .
$$
Now, if $\|X\|_{2p, w} < \infty$, we get from Theorem \ref{order1} that
\[
 \left \| W_p^p(\mu_n, \mu)  \right \|_1 \ll  \left\{ \begin{aligned}
 & 
 \frac{\|X\|_{2p,w}^{p} \log n}{\sqrt n} 
 & \text{ if } & p > d/2 \\
 & \frac{\|X\|_{2p,w}^{p} (\log n)^{2}}{\sqrt n} & \text{ if } & p = d/2 \\
 &  \frac{\|X\|_{2p,w}^{p} }{n^{p/d}} 
 & \text{ if } & p  \in [1, d/2) \, .
    \end{aligned}
    \right.
\]
Finally, if $\int_0^{\infty} t^{p-1} \sqrt{H (t) }  \,  dt < \infty$, the rates in the cases $p>d/2$ and $p=d/2$  can be slightly improved
 (taking $q=2$ and $M=\infty$ in Theorem \ref{order1}); this can be directly deduced from 
Theorem \ref{order2} below. 

Note that all those bounds are consistent with that given in Theorem 1
of \cite{FG}, and slightly more precise in terms of the moment conditions. Hence, the
discussion on the optimality of the rates in \cite{FG} is also valid for our Theorem \ref{order1} (see Remark \ref{momentOpt} below). 
For $p<d/2$ and $\|X\|_{q}< \infty $ for some $q>dp/(d-p)$, it follows from Theorem 2(ii) in \cite{DSS} that 
$\liminf_{n \rightarrow \infty} n^{p/d} \left \| W_p^p(\mu_n, \mu)  \right \|_1>0$ if $\mu$ has a non degenerate absolutely continuous part
with respect to the Lebesgue measure, and that $\limsup_{n \rightarrow \infty} n^{p/d} \left \| W_p^p(\mu_n, \mu)  \right \|_1=0$ if $\mu$ is singular. 
Still for $p<d/2$, we refer to the paper \cite{BW}, which shows that, for compactly supported singular measures, the rates of convergence  of 
 $\| W_p^p(\mu_n, \mu)   \|_1$
can be much faster than  $n^{-p/d}$. 
\end{Remark}

\begin{Theorem}\label{order2}
If $\int_0^{\infty} t^{p-1} \sqrt{H (t) }  \,  dt < \infty$, then
\[
 \left \| W_p^p(\mu_n, \mu)  \right \|_2^2 \ll  \left\{ \begin{aligned}
 & \frac{1}{n}  
 \left (\int_0^{\infty} t^{p-1} \sqrt{ H(t)} \,     dt  \right )^2 
 & \text{ if } & p > d/2 \\
 &\frac{(\log n)^2}{n}  
 \left (\int_0^{\infty} t^{p-1}\sqrt{ H(t)} \,    dt  \right )^2 & \text{ if } & p = d/2 \\
 &   \frac{1}{n^{2p/d} }
 \left (\int_0^{\infty} t^{p-1}\sqrt{ H(t)} \,     dt  \right )^2 
 & \text{ if } & p  \in [1, d/2)
    \end{aligned}
    \right.
\]
\end{Theorem}

\begin{Remark}\label{momentOpt}
 According to the discussion after Theorem 1 in \cite{FG}, if $p \neq d/2$, one can always find some measure $\mu$ 
for which the rate of Theorem \ref{order2} is reached (see example (a) and (b)   in \cite{FG} for $ p > d/2$ and example (c) 
 in \cite{FG} for $ p < d/2$). 
 
Note also that, for ${\mathbb E}(W_1(\mu_n, \mu))$ instead of $\|W_1(\mu_n, \mu)\|_2$,  the bounds of Theorem \ref{order2}  can be obtained   from the general bound given in Theorem 3.8
of \cite{LG}, under the condition $\int_0^{\infty} t^{p-1} \sqrt{H (t) }  \,  dt < \infty$ (taking a ball of radius $r=H^{-1}(\alpha)$ to bound up the 
term $\tau_n^\alpha$ in \cite{LG}, and noting that $\int_0^{\infty} t^{p-1} \sqrt{H (t) }  \,  dt < \infty$ is equivalent to $\int_0^1 (H^{-1}(\alpha))^p \alpha^{-1/2}\,  d\alpha < \infty$).
 
In the case $d=1, p=1$, del Barrio et al. \cite{dBGM} proved that $ \sqrt n W_1(\mu_n, \mu)$
is stochastically bounded if and only if $ \int_0^{\infty}  \sqrt{H (t) }  \,  dt < \infty$ (see their Theorem 2.1(b)), which is consistent with the first inequality 
of Theorem \ref{order2}. 
For $d=1, p>1$,  we refer to the paper by Bobkov and Ledoux \cite{BL} for some conditions on $\mu$ ensuring faster rates of convergence.
Finally, when $p=1, d=2$ and $\mu$ is the uniform measure over $[0,1]^2$, Ajtai et al. \cite{AKT} proved that 
 ${\mathbb E}(W_1(\mu_n, \mu))$ is exactly of order $(\log n/n)^{1/2}$, while we get a rate of order $ \log n / \sqrt n$, which is therefore suboptimal in that 
particular case.  
For other  discussions about the rates, see for instance \cite{K}, Sections 2.3 and 2.4.
\end{Remark}

\subsection{von Bahr-Esseen type inequalities}

In this subsection, we shall prove some moment inequalities in the spirit of von Bahr and Esseen \cite{vBE}. Recall that, for partial sums
$S_n=Y_1+\cdots + Y_n$ of i.i.d real-valued random variables such that $\|Y_1\|_r <\infty$  for some $r \in [1,2]$ and ${\mathbb E}(Y_1)=0$, the inequality 
of von Bahr and Esseen reads as follows:
\begin{equation}\label{Ivbe}
 \left \| \frac{S_n}{n} \right \|_r^r  \leq 
  \frac{2\|Y_1\|_{r}^{r}}{n^{r-1}}  \, .
\end{equation}

\medskip

In the case  case where $r \in (1,2)$, we prove the following result. 

\begin{Theorem}\label{VBE}
If $\|X\|_{rp} < \infty$ for some $r \in (1,2)$, then
\[
 \left \| W_p^p(\mu_n, \mu)  \right \|_r^r  \ll  \left\{ \begin{aligned}
 & \frac{\|X\|_{rp}^{rp}}{n^{r-1}}  & \text{ if } & p > d(r-1)/r \\
 &   \frac{\|X\|_{rp}^{rp}}{n^{rp/d}}   & \text{ if } & p  \in [1, d(r-1)/r) 
    \end{aligned}
    \right.
\]
\end{Theorem}

\begin{Remark}
For $d=1$, the first inequality of Theorem \ref{VBE} has been proved in \cite{DM}. 
Our proof does not allow to deal with the case where $p=d(r-1)/r$. However, in that case, it is easy to see 
that 
$$
 \left \| W_p^p(\mu_n, \mu)  \right \|_r^r  \leq  \frac{(\log n)^r}{n^{r-1}}
 \left (\int_0^{\infty} t^{p-1} (H(t))^{1/r} \,    dt  \right )^r
$$
(same proof as the second inequality of Theorem \ref{order2}). 
For $p=1$ and $d<r/(r-1)$, using the dual expression of $W_1(\mu_n, \mu)$ (see \eqref{dual}), we get the upper bound
\begin{equation}\label{unif}
  \left \| \sup_{f \in \Lambda_1} \left | \frac 1 n  \sum_{k=1}^n \left  ( f(X_k) - \mu(f)  \right ) \right |  \, \right \|_r^r \ll \frac{\|X\|_{r}^{r}}{n^{r-1}} \, ,
\end{equation}
where $\Lambda_1$ is the the set of functions $f$ such that $|f(x)-f(y)| \leq |x-y|_2$.
Note that \eqref{unif} may be seen as a uniform version of the inequality \eqref{Ivbe} over the class $\Lambda_1$. 
\end{Remark} 

\subsection{Rosenthal type inequalities}
In this subsection, we shall prove some moment inequalities in the spirit of Rosenthal \cite{Ro}. Recall that, for partial sums
$S_n=Y_1+\cdots + Y_n$ of i.i.d real-valued random variables such that $\|Y_1\|_r <\infty$  for some $r \geq 2$ and ${\mathbb E}(Y_1)=0$, the inequality 
of Rosenthal reads as follows: there exists two positive constants $c_1(r)$ and $c_2(r)$ such that
\begin{equation*}\label{Iros}
 \left \| \frac{S_n}{n} \right \|_r^r  \leq 
   c_1(r) \frac{\|Y_1\|_{2}^{r}}{n^{r/2}}     + c_2(r)\frac{\|Y_1\|_{r}^{r}}{n^{r-1}}\, .
\end{equation*}
We refer to Pinelis \cite{P}  for the expression of the possible constants $c_1(r)$ and $c_2(r)$. 

\medskip

In  the case where $r>2$, we prove the following result.

\begin{Theorem}\label{Rosine}
If $\|X\|_{rp} < \infty$ for some $r >2$, then
\[
 \Vert  W_p^p (\mu_n, \mu) \Vert_r^r \ll    \left\{ \begin{aligned}
 &   \frac{1}{n^{r/2}}   \left (\int_0^{\infty} t^{p-1} \sqrt{H  (t)}  \,   dt  \right)^r + \frac{\|X\|_{rp}^{rp}}{n^{r-1}}    & \text{ if } & p  > d(r-1)/r   \\
 &\frac{1}{n^{r/2}}   \left (\int_0^{\infty} t^{p-1} \sqrt{H  (t)}   \,  dt  \right )^r + \frac{n^{\gamma }} {n^{pr/d} } \|X\|_{rp}^{rp} & \text{ if } & d/2 < p  \leq d(r-1)/r  \\
 & \frac{(\log n)^r}{n^{r/2}}  \left(\int_0^{\infty} t^{d/2 -1} \sqrt{H(t)} dt \right )^r  + \frac{(\log n)^2}{n^{r/2}}   \|X\|_{rp}^{rp} & \text{ if } & p = d/2 \\
 & \frac{\|X\|_{rp}^{rp}}{n^{rp/d}} & \text{ if } & 
 p  \in [1, d/2) 
    \end{aligned}
    \right.
\]
where, for the second inequality,  $\gamma$ can be taken as $\gamma =  \frac{\varepsilon (2p-d)}{ d   (r-2 + \varepsilon)}$ for any $\varepsilon >0$ (and the constants 
implicitely involved in the inequality depend on $\varepsilon$). 
\end{Theorem}

\begin{Remark}
For $d=1$, the first inequality of Theorem \ref{Rosine} has been proved in \cite{DM}. 
As a consequence of the two first inequalities of Theorem \ref{Rosine}, we obtain that, if $p>d/2$,
$$
\limsup_{n \rightarrow \infty} \sqrt n  \Vert  W_p^p (\mu_n, \mu) \Vert_r \ll \int_0^{\infty} t^{p-1} \sqrt{H  (t)}  \,   dt   \, .
$$
As a consequence of the third inequality of Theorem \ref{Rosine}, we obtain that,
if $p=d/2$,
$$
\limsup_{n \rightarrow \infty} \frac{\sqrt n}{\log n}  \Vert  W_p^p (\mu_n, \mu) \Vert_r \ll \int_0^{\infty} t^{p-1} \sqrt{H  (t)}  \,   dt   \, .
$$
Note also that,  according to the discussion after Theorem 1 in \cite{FG}, if $p \neq d/2$, one can always find some measure $\mu$ 
for which the rates of Theorem \ref{Rosine} are reached (see example (a)  in \cite{FG} for $ p > d/2$ and example (c) 
 in \cite{FG} for $ p < d/2$).

\end{Remark}

\section{Proofs}\label{Proofs}

\setcounter{equation}{0}

The starting point of the proofs is Lemmas 5 and  6 in \cite{FG}, which we recall below. 

For $\ell \geq 0$,  let ${\mathcal P}_{\ell}$ be the natural partition of $(-1,1]^d$ into $2^{d \ell}$ translations
of $(-2^{-\ell}, 2^{-\ell}]^d$.  Let also $B_0=(-1,1]^d$ and for any integer $m\geq 1$, 
$B_m= (-2^m, 2^m]^d \,  \backslash \,  (-2^{m-1}, 2^{m-1}]^d$.
For a set $F \subset {\mathbb R}^d$ and $a>0$, we use the standard notation $aF=\{ax : x \in F\}$.  For a probability measure $\nu$ on ${\mathbb R}^d$ and $m \geq 0$, let ${\mathcal R}_{B_m} \nu$ be the probability measure on $(-1,1]^d$ defined as the image of $\nu |_{B_m} / \nu (B_m)$ by the map $x \mapsto x/2^m$.  For two probability measures $\mu$ and $\nu$ on ${\mathbb R}^d$, by Lemma 5 in \cite{FG}, there exists a positive constant $\kappa_{p,d}$ depending only on $p$ and $d$ such that 
\beq \label{lma5FG}
W_p^p ( \mu, \nu) \leq \kappa_{p,d} {\mathcal D}_p (\mu,\nu) \, , 
\eeq
where 
\beq \label{defclaDp}
{\mathcal D}_p (\mu,\nu) := \sum_{m \geq 0} 2^{pm} |  \mu (B_m) - \nu (B_m) | +  \sum_{m \geq 0} 2^{pm} (  \mu (B_m) \wedge \nu (B_m) ) {\mathcal D}_p ({\mathcal R}_{B_m}\mu,{\mathcal R}_{B_m} \nu)  \, ,
\eeq
with
\beq \label{defclaDpR}
{\mathcal D}_p ({\mathcal R}_{B_m}\mu,{\mathcal R}_{B_m} \nu)  = \frac{2^p -1}{2} \sum_{\ell \geq 1} 2^{-p \ell} \sum_{F \in {\mathcal P}_{\ell}} 
\left |  \frac{\mu( 2^m F \cap B_m ) }{\mu (B_m) } -  \frac{\nu ( 2^m F \cap B_m )}{\nu (B_m) }\right | \, .
\eeq
In addition, by Lemma 6 in \cite{FG},
\[
{\mathcal D}_p (\mu,\nu) \leq \left ( \frac{3}{2} \vee  \frac{2^{p} -1}{2}  \right ) \Delta_p (\mu, \nu ) 
\]
where
\[
\Delta_p (\mu, \nu ) = \sum_{m \geq 0 } 2^{pm} \sum_{\ell \geq 0} 2^{-p \ell} \sum_{ F \in {\mathcal P}_{\ell}} | \mu ( 2^m F \cap B_m ) - \nu ( 2^m F \cap B_m )  | \, . 
\]
From the considerations above, there exists a constant $C$ depending only on $p$ and $d$ such that 
\begin{equation}\label{depart}
W_p^p (\mu_k, \mu )  \leq  C \Delta_p (\mu_k, \mu )  \, ,
\end{equation}
where $\mu_k = \frac{1}{k} \sum_{i=1}^k \delta_{X_i}$. 
This inequality may be seen as an extension to the case $d>1$ of \`Ebralidze's inequality \cite{E}, which we used  in \cite{DM} to obtain moment bounds for $W_p^p (\mu_n, \mu )$ when 
$d=1$. 

\medskip

As in \cite{DM}  we shall use truncation arguments. 
For a positive real  $M$, let ${\mathcal C}_M =[-M,M]^d$, 
\[
A_{p,M} (\mu_k, \mu ) = \sum_{m \geq 0 } 2^{pm} \sum_{\ell \geq 0} 2^{-p \ell} \sum_{ F \in {\mathcal P}_{\ell}} | \mu_k ( 2^m F \cap B_m \cap  {\mathcal C}_M) - \mu ( 2^m F \cap B_m \cap  {\mathcal C}_M )  | 
\]
and
\[
B_{p,M} (\mu_k, \mu ) = \sum_{m\geq 0 } 2^{pm} \sum_{\ell \geq 0} 2^{-p \ell} \sum_{ F \in {\mathcal P}_{\ell}} | \mu_k ( 2^m F \cap B_m \cap  {\mathcal C}^c_M) - \mu ( 2^m F \cap B_m \cap  {\mathcal C}^c_M )  |  \, .
\]
With these notations, 
it follows  that
\begin{equation}\label{departbis}
\Delta_p (\mu_k, \mu )  \leq A_{p,M} (\mu_k, \mu )  + B_{p,M} (\mu_k, \mu ) \, .
\end{equation}

For the proofs, we shall follow the order of the theorems, except for Theorem \ref{VBE} whose proof comes naturally after those of Theorems \ref{weakVBE} and \ref{F-Nag}.

\subsection{Proof of Theorem \ref{weakVBE}}\label{sec:WVBE}

Let $M>0$ and $x >0$. Starting from \eqref{depart} and \eqref{departbis}, we get that
\begin{multline}\label{start}
{\mathbb P}\left ( \max_{1 \leq k \leq n}kW_p^p(\mu_k, \mu) >nx \right )  \leq 
{\mathbb P}\left ( \max_{1 \leq k \leq n}kA_{p,M}(\mu_k, \mu) >(nx/2C)\right ) \\ + 
{\mathbb P}\left ( \max_{1 \leq k \leq n}kB_{p,M}(\mu_k, \mu) >(nx /2C)\right ) \, .
\end{multline}

Let $y=x/2C$. By Markov's inequality at order  $q \in  (r, 2)$ and  $s \in [1,r)$,
\begin{align}
{\mathbb P}\left ( \max_{1 \leq k \leq n}kA_{p,M}(\mu_k, \mu) >ny\right )  & \leq  \frac{ \left \| \max_{1 \leq k \leq n}kA_{p,M}(\mu_k, \mu)\right \|_q^q}{n^q y^q}  \, , \label{A_p}\\
{\mathbb P}\left ( \max_{1 \leq k \leq n}kB_{p,M}(\mu_k, \mu) >ny\right ) & \leq \frac{\left \|  \max_{1 \leq k \leq n}kB_{p,M}(\mu_k, \mu)\right \|_s^s}{n^s y^s} \label{B_p} \, .
\end{align}

To deal with  \eqref{A_p},  we first note that 
\begin{multline} \label{decA-16-dec}
\Big \|  \max_{1 \leq k \leq n}kA_{p,M}(\mu_k, \mu) \Big  \|_q  \\ \ll 
\sum_{m\geq 0 } 2^{pm} \sum_{\ell \geq 0} 2^{-p \ell} \Big  \|  \max_{1 \leq k \leq n}\sum_{ F \in {\mathcal P}_{\ell}} | k\mu_k ( 2^m F \cap B_m \cap  {\mathcal C}_M) - k\mu ( 2^m F \cap B_m \cap  {\mathcal C}_M )  |  \Big  \|_q \, .
\end{multline}
Now, clearly
\begin{multline}\label{step1}
\Big  \|  \max_{1 \leq k \leq n}\sum_{ F \in {\mathcal P}_{\ell}} | k\mu_k (  2^m F \cap B_m \cap  {\mathcal C}_M) - k\mu ( 2^m F \cap B_m \cap  {\mathcal C}_M )  |  \Big  \|_q \\
\leq \Big  \|  \max_{1 \leq k \leq n} \left ( k\mu_k (   B_m \cap  {\mathcal C}_M) +  k\mu (  B_m \cap  {\mathcal C}_M )  \right )  \Big  \|_q 
\leq 2 n \left( \mu (  B_m \cap  {\mathcal C}_M ) \right )^{1/q} \, .
\end{multline}
On the other hand,  by the (maximal version of) von Bahr-Essen inequality (see \cite{vBE}),
$$
\Big  \| \max_{1 \leq k \leq n}| k\mu_k ( 2^m F \cap B_m \cap  {\mathcal C}_M) - k\mu ( 2^m F \cap B_m \cap  {\mathcal C}_M )  | \Big  \|_q^q
\ll n  \mu ( 2^m F \cap B_m \cap  {\mathcal C}_M ) \, ,
$$
so that, by using H\"older's inequality and the fact that $| {\mathcal P}_{\ell}|=2^{\ell d}$, 
\begin{equation}\label{step2}
\sum_{ F \in {\mathcal P}_{\ell}} \Big  \| \max_{1 \leq k \leq n} | k\mu_k (  2^m F \cap B_m \cap  {\mathcal C}_M) - k\mu ( 2^m F \cap B_m \cap  {\mathcal C}_M )  | \Big  \|_q 
\ll  2^{\ell d (q-1)/q} n^{1/q} \left( \mu (  B_m \cap  {\mathcal C}_M ) \right )^{1/q}  \, .
\end{equation}
Combining \eqref{A_p}, \eqref{step1} and \eqref{step2}, we obtain that 
\begin{equation}\label{Abound}
{\mathbb P}\left ( \max_{1 \leq k \leq n}kA_{p,M}(\mu_k, \mu) >ny\right ) \ll 
\frac{1}{y^q} \left ( \sum_{m\geq 0 } 2^{pm} \left( \mu (  B_m \cap  {\mathcal C}_M ) \right )^{1/q} \sum_{\ell \geq 0}   \frac {1}{2^{p \ell}}    \min \left (1, n^{-(q-1)/q}2^{\ell d (q-1)/q} \right) \right)^q .
\end{equation}

In the same way, for the term \eqref{B_p}, we obtain the upper bound 
\begin{equation}\label{Bbound}
{\mathbb P}\left ( \max_{1 \leq k \leq n}kB_{p,M}(\mu_k, \mu) >ny\right ) \ll 
\frac{1}{y^s} \left ( \sum_{m\geq 0 } 2^{pm} \left( \mu (  B_m \cap  {\mathcal C}^c_M ) \right )^{1/s} \sum_{\ell \geq 0}  \frac {1}{2^{p \ell}}  \min \left (1, n^{-(s-1)/s}2^{\ell d (s-1)/s} \right) \right)^s .
\end{equation}

From \eqref{Abound} and \eqref{Bbound}, we see that three cases arise:

\medskip

\noindent $\bullet$  If $p>d(r-1)/r$, then, taking  $q>r$ such that $p>d(q-1)/q$ and $s=1$, we get the upper bounds
\begin{multline}\label{Abound1}
{\mathbb P}\left ( \max_{1 \leq k \leq n}kA_{p,M}(\mu_k, \mu) >ny\right ) \\ \ll 
\frac{1}{n^{q-1}y^q} \left ( \sum_{m\geq 0 } 2^{pm}   \left( \mu (  B_m \cap  {\mathcal C}_M ) \right )^{1/q}  \right )^q  \ll \frac{1}{n^{q-1}y^q} 
 \left( \int_0^\infty t^{p-1} (H(t))^{1/q} {\bf 1}_{t \leq M} dt
 \right )^q \, , 
\end{multline}
and 
\begin{equation}\label{Bbound1}
{\mathbb P}\left ( \max_{1 \leq k \leq n}kB_{p,M}(\mu_k, \mu) >ny\right ) \ll \frac{1}{y} \int_0^\infty t^{p-1} H(t) {\bf 1}_{t >M} dt \, .
\end{equation}

Using that $H(t) \leq \|X\|_{rp, w}^{rp} t^{-r p} $ for $r \in (1,2)$, we infer from  \eqref{start},  \eqref{Abound1} and \eqref{Bbound1} that
$$
{\mathbb P}\left ( \max_{1 \leq k \leq n}kW_p^p(\mu_k, \mu) >nx \right )  \ll  \|X\|_{rp, w}^{rp} \left (  \frac{1}{x M^{p(r-1)}} + \frac{M^{p(q-r)}}{n^{q-1} x^q}\right )  \, .
$$
Taking $M=(nx)^{1/p}$, we obtain the desired result when $p>d(r-1)/r$. 

\smallskip

\noindent $\bullet$  If $p=d(r-1)/r$, then, taking $q=r$, we get the upper bound
\begin{multline}\label{Abound2}
{\mathbb P}\left ( \max_{1 \leq k \leq n}kA_{p,M}(\mu_k, \mu) >ny\right ) \\ \ll 
\frac{(\log(n))^r}{n^{r-1}y^r} \left ( \sum_{m\geq 0 } 2^{pm}   \left( \mu (  B_m \cap  {\mathcal C}_M ) \right )^{1/r}  \right )^r  \ll \frac{(\log n)^r}{n^{r-1}y^r}
 \left( \int_0^\infty t^{p-1} (H(t))^{1/r} {\bf 1}_{t \leq M} dt
 \right )^r \, .
\end{multline}
Using that $H(t) \leq \|X\|_{rp, w}^{rp} t^{-r p} $ for $r \in (1,2)$, we infer from \eqref{start}, \eqref{Bbound1} and  \eqref{Abound2}  that
$$
{\mathbb P}\left ( \max_{1 \leq k \leq n}kW_p^p(\mu_k, \mu) >nx \right )  \ll  \|X\|_{rp, w}^{rp} \left (  \frac{1}{x M^{p(r-1)}} + \frac{(\log n)^r}{n^{r-1} x^r}
\left( 1+ \log_+\left(  \frac{M}{ \|X\|_{rp,w} }\right) \right)^r \right )  \, .
$$
Taking $M=(nx)^{1/p}$, we obtain the desired result when $p=d(r-1)/r$. 

\smallskip

\noindent $\bullet$  If $p<d(r-1)/r$, then, taking $q>r$ and $s \in (1,r)$ such that $p<d(s-1)/s$, we get the upper bounds
\begin{multline}\label{Abound3}
{\mathbb P}\left ( \max_{1 \leq k \leq n}kA_{p,M}(\mu_k, \mu) >ny\right ) \\  \ll 
\frac{1}{n^{qp/d}y^q} \left ( \sum_{m\geq 0 } 2^{pm}   \left( \mu (  B_m \cap  {\mathcal C}_M ) \right )^{1/q}  \right )^q  \ll \frac{1}{n^{qp/d}y^q} 
 \left( \int_0^\infty t^{p-1} (H(t))^{1/q} {\bf 1}_{t \leq M} dt
 \right )^q \, ,
\end{multline}
and 
\begin{multline}\label{Bbound3}
{\mathbb P}\left ( \max_{1 \leq k \leq n}kB_{p,M}(\mu_k, \mu) >ny\right ) \\  \ll 
\frac{1}{n^{sp/d}y^s} \left ( \sum_{m\geq 0 } 2^{pm}   \left( \mu (  B_m \cap  {\mathcal C}^c_M ) \right )^{1/s}  \right )^s  \ll \frac{1}{n^{sp/d}y^s} 
 \left( \int_0^\infty t^{p-1} (H(t))^{1/s} {\bf 1}_{t >M} dt
 \right )^s \, ,
\end{multline}

Using that $H(t) \leq \|X\|_{rp, w}^{rp} t^{-r p} $ for $r \in (1,2)$, we infer from \eqref{start}, \eqref{Abound3}  and  \eqref{Bbound3} that
$$
{\mathbb P}\left ( \max_{1 \leq k \leq n}kW_p^p(\mu_k, \mu) >nx \right )  \ll  \|X\|_{rp, w}^{rp} \left (  \frac{1}{n^{sp/d} x^s M^{p(r-s)}} + \frac{M^{p(q-r)}}{n^{qp/d} x^q}
 \right )  \, .
$$
Taking $M=n^{1/d}x^{1/p}$, we obtain the desired result when $p<d(r-1)/r$. 

\subsection{Proof of Theorem \ref{F-Nag}}
Let $r>2$. Note first that, by homogeneity, the general inequality may be deduced from the case where $\|X\|_{rp, w}=1$ by considering the variables 
$X_i/\|X\|_{rp, w}$. Hence, from now, we shall assume that $\|X\|_{rp, w}=1$. 

According to the beginning of the proof of Theorem 2 in \cite{FG}, we get that 
\begin{equation}\label{fgu1}
W_p^p(\mu_n, \mu) \leq C \sum_{m \geq 0}  2^{pm} \left | \mu_n(B_m)- \mu(B_m) \right | + C V_n^p \, ,
\end{equation}
for some positive constant $C=C_{p,d}$, where the random variable $V_n^p$ is such that 
\begin{equation}\label{fgu2}
{\mathbb P}( V_n^p \geq x/(2C)) \leq a(n,x) \, .
\end{equation}
Consequently, it remains to bound up the quantity
$$
{\mathbb P} \left ( \sum_{m \geq 0}  2^{pm} \left | \mu_n(B_m)- \mu(B_m) \right | \geq  x/(2C) \right )  \, .
$$

For a positive real  $M$, let ${\mathcal C}_M =[-M,M]^d$, 
\[
A^*_{p,M} (\mu_k, \mu ) = \sum_{m \geq 0 } 2^{pm}  | \mu_k ( B_m \cap  {\mathcal C}_M) - \mu ( B_m \cap  {\mathcal C}_M )  | 
\]
and
\[
B^*_{p,M} (\mu_k, \mu ) = \sum_{m\geq 0 } 2^{pm}  | \mu_k ( B_m \cap  {\mathcal C}^c_M) - \mu ( B_m \cap  {\mathcal C}^c_M )  |  \, .
\]
With these notations,
\begin{multline}\label{start2}
{\mathbb P} \left ( \sum_{m \geq 0}  2^{pm} \left | \mu_n(B_m)- \mu(B_m) \right | \geq  x/(2C) \right )  \leq 
{\mathbb P}\left ( A^*_{p,M}(\mu_n, \mu) >x/(4C)\right ) \\ + 
{\mathbb P}\left ( B^*_{p,M}(\mu_n, \mu) >x /(4C)\right ) \, .
\end{multline}
Let $y=x/4C$. By Markov's inequality at order  $q>2$ and 1,
\begin{align}
{\mathbb P}\left ( A^*_{p,M}(\mu_n, \mu) >y\right )  & \leq  \frac{ \left \| A^*_{p,M}(\mu_n, \mu)\right \|_q^q}{y^q}  \, , \label{A_p2}\\
{\mathbb P}\left ( B^*_{p,M}(\mu_n, \mu) >y\right ) & \leq \frac{\left \|  B^*_{p,M}(\mu_n, \mu)\right \|_1}{y} \label{B_p2} \, .
\end{align}
 Applying Rosenthal's inequality, we get
\begin{multline*}
 \Vert A^*_{p,M} (\mu_n, \mu ) \Vert_q \ll \frac{1}{\sqrt n}
  \sum_{m \geq 0 } 2^{pm}  \mu^{1/2} ( B_m  \cap  {\mathcal C}_M )   +  
   \frac{1}{n^{(q-1)/q}}\sum_{m \geq 0 } 2^{pm}  \mu^{1/q} ( B_m  \cap  {\mathcal C}_M )   \\
 \ll \frac{1}{\sqrt n} \int_0^{\infty} t^{p-1} \sqrt{H (t)} {\bf 1}_{t \leq M} dt + 
 \frac{1}{n^{(q-1)/q}}  \int_0^{\infty} t^{p-1} H^{1/q} (t) {\bf 1}_{t \leq M} dt \, .
 \end{multline*}
Choosing $q >r$, 
it follows that
\begin{align}
 \Vert A^*_{p,M} (\mu_n, \mu ) \Vert_q 
 & \ll \frac{1}{\sqrt n}  \int_0^{\infty} t^{p-1} \sqrt{H (t)}  dt + \frac{M ^{ p (q-r)/q }}{n^{(q-1)/q}}   \nonumber
  \left ( \sup_{t >0}t^{rp} H(t)  \right )^{1/q} \\
 & \ll \frac{1}{\sqrt n}  \int_0^{\infty} t^{p-1} \sqrt{H (t)}  dt + \frac{M ^{ p (q-r)/q }}{n^{(q-1)/q}}   \, , \label{ouf2}
\end{align}
the last inequality being true since we assumed that $\sup_{t >0}t^{rp} H(t)=1$. 

On another hand, 
\begin{align}
 \Vert B^*_{p,M} (\mu_n, \mu ) \Vert_1 &\ll    \sum_{n \geq 0 } 2^{pn}  \mu ( B_n  \cap  {\mathcal C}^c_M ) \nonumber
     \\
 &\ll  \int_0^{\infty} t^{p-1} H (t) {\bf 1}_{t > M} dt  \ll M^{p(1-r)}
  \sup_{t >0}t^{rp} H(t)    \ll  M^{p(1-r)}\, .  \label{ouf}
 \end{align}
 Gathering \eqref{start2}  - \eqref{ouf}, we get that for any $q >r$,
\begin{equation}\label{ouff}
 {\mathbb P}\left ( \sum_{m \geq 0} 2^{pm} \left | \mu_n (B_m) - \mu (B_m)  \right |  > x/(2C) \right )  
  \ll  \frac{1}{x^{q}n^{q/2}}  \left ( \int_0^{\infty} t^{p-1} \sqrt{H (t)}  dt  \right )^q  + \frac{ M^{ p(q   - r ) }}{x^{q}n^{q-1}}   + \frac{M^{p(1-r)}}{x}     \, ,
 \end{equation}
 Hence choosing $M = n^{1/p}x^{1/p}$, we infer from \eqref{fgu1}, \eqref{fgu2} and \eqref{ouff} that for  any  $q >r$,
\begin{equation} \label{FG1}
 {\mathbb P} \left (  W^p_{p} (\mu_n, \mu ) >  x \right )   \ll a(n,x)   +  \frac{1}{x^r n^{r -1}}  +  \frac{1}{x^q n^{q/2}}  \left ( \int_0^{\infty} t^{p-1} \sqrt{H (t)} dt \right )^q  \, ,
 \end{equation}
 which is the desired inequality when $\sup_{t >0}t^{rp} H(t)=\|X\|_{rp, w}^{rp}=1$.

\subsection{Proof of Theorem \ref{VBE}}

We start from the elementary equality
\begin{equation}\label{start0}
\frac{1}{n^r}\left \| \max_{1 \leq k \leq n}kW_p^p(\mu_k, \mu) \right \|_r^r= r\int_0^\infty x^{r-1} {\mathbb P}\left ( \max_{1 \leq k \leq n}kW_p^p(\mu_k, \mu) >nx \right ) dx \, .
\end{equation}
Then, we  use the upper bounds \eqref{start}, \eqref{Abound} and \eqref{Bbound}. We consider  two cases:

\medskip

\noindent $\bullet$  If $p>d(r-1)/r$,  let $q>r$ such that $p>d(q-1)/q$, and let $M=(nx)^{1/p}$. From \eqref{start}, \eqref{Abound1} and \eqref{Bbound1} we get the upper bound
\begin{multline}\label{M1}
\int x^{r-1} {\mathbb P}\left ( \max_{1 \leq k \leq n}kW_p^p(\mu_k, \mu) >nx \right ) dx  \\ \ll  \int_0^\infty x^{r-2} \left(\int_0^\infty t^{p-1} H(t) {\bf 1}_{t >M} dt \right) dx + \frac{1}{n^{q-1}}  \int_0^\infty x^{r-1-q}
 \left( \int_0^\infty t^{p-1} (H(t))^{1/q} {\bf 1}_{t \leq M}  dt \right )^q dx \, .
\end{multline}
Note that 
\begin{equation}\label{M2}
\int_0^\infty x^{r-2} \left(\int_0^\infty t^{p-1} H(t) {\bf 1}_{t >(nx)^{1/p}} dt \right) dx \ll  \frac{1}{n^{r-1}}\int_0^\infty t^{rp-1} H(t) dt \ll \frac{\|X\|_{rp}^{rp}}{n^{r-1}}
\, .
\end{equation}
Let $\beta < (q-1)/q$. Applying H\"older's inequality, we obtain
\begin{multline}\label{Heq}
\frac {1}{ n^{q-1} }\int_{0}^{\infty} x^{r-1-q }  \left (  \int_0^{\infty} t^{p-1} (H (t))^{1/q} {\bf 1}_{t \leq (nx)^{1/p}}  dt  \right )^q  dx \\
  \ll  \frac {1}{n^{q-1}}  \int_{0}^{\infty} x^{r-1-q} (nx)^{(q-1 - q\beta) /p}  \left ( \int_0^{\infty} t^{q(p-1 + \beta ) } H (t) {\bf 1}_{t \leq (nx)^{1/p}}  dt  \right )  dx  \\
  \ll  \frac{ n^{(q-1 - q\beta) /p}}{n^{q-1}}\int_0^{\infty} t^{q(p-1 + \beta ) } H (t) \left ( \int_0^{\infty } x^{r-1-q +(q-1 - q\beta) /p }     {\bf 1}_{ x \geq t^p /n }   dx  \right )dt    \, .
\end{multline}
Taking $\beta$ close enough to $(q-1)/q$ in such a way that $q+1 -r-(q-1 - q\beta) /p >1$, we get
\begin{equation}\label{M3}
\frac {1}{ n^{q-1} }\int_{0}^{\infty} x^{r-1-q }  \left (  \int_0^{\infty} t^{p-1} (H (t))^{1/q} {\bf 1}_{t \leq (nx)^{1/p}}  dt  \right )^q  dx  \ll \frac{\|X\|_{rp}^{rp}}{n^{r-1}}
\, .
\end{equation}
Gathering \eqref{M1}, \eqref{M2} and \eqref{M3}, we obtain the desired result.


\smallskip

\noindent $\bullet$  If $p<d(r-1)/r$,  let  $q>r$,  $s  \in (1,r) $ such that $p<d(s-1)/s$, and let $M=n^{1/d}x^{1/p}$. From \eqref{start},  \eqref{Abound3} \eqref{Bbound3} we get the upper bound
\begin{multline}\label{M1bis}
\int x^{r-1} {\mathbb P}\left ( \max_{1 \leq k \leq n}kW_p^p(\mu_k, \mu) >nx \right ) dx   \ll  \frac{1}{n^{sp/d}}  \int_0^\infty x^{r-1-s}
 \left( \int_0^\infty t^{p-1} (H(t))^{1/s} {\bf 1}_{t > M}  dt \right )^s dx 
 \\ + 
\frac{1}{n^{qp/d}}  \int_0^\infty x^{r-1-q}
 \left( \int_0^\infty t^{p-1} (H(t))^{1/q} {\bf 1}_{t \leq M}  dt \right )^q dx \, .
\end{multline}
Proceeding exactly as for \eqref{Heq}-\eqref{M3}, with the choice $M=n^{1/d}x^{1/p}$, we get
\begin{equation}\label{M4}
\frac {1}{ n^{pq/d} }\int_{0}^{\infty} x^{r-1-q }  \left (  \int_0^{\infty} t^{p-1} (H (t))^{1/q} {\bf 1}_{t \leq n^{1/d}x^{1/p}}  dt  \right )^q  dx  \ll 
\frac{1}{n^{pr/d}}\int_0^\infty t^{rp-1} H(t) dt \ll \frac{\|X\|_{rp}^{rp}}{n^{pr/d}}
\, .
\end{equation}
In the same way, we get 
\begin{equation}\label{M4bis}
\frac {1}{ n^{ps/d} }\int_{0}^{\infty} x^{r-1-s }  \left (  \int_0^{\infty} t^{p-1} (H (t))^{1/s} {\bf 1}_{t >  n^{1/d}x^{1/p}}  dt  \right )^s  dx  \ll 
\frac{1}{n^{pr/d}}\int_0^\infty t^{rp-1} H(t) dt \ll \frac{\|X\|_{rp}^{rp}}{n^{pr/d}}
\, .
\end{equation}
Gathering \eqref{M1bis}, \eqref{M4} and \eqref{M4bis}, we obtain the desired result.

\subsection{Proof of Theorem \ref{BK1}} \label{proofbk}

Let $r \in (1,2)$. We  start from the upper bounds \eqref{start}, \eqref{Abound} and \eqref{Bbound}. 

\medskip

\noindent $\bullet$  If $p>d(r-1)/r$, let $q \in (r,2]$ such that $p>d(q-1)/q$ and  let $M=n^{\alpha/p}$. From \eqref{start}, \eqref{Abound1}  and  \eqref{Bbound1} we get the upper bound
$$
{\mathbb P}\left ( \max_{1 \leq k \leq n} k W_p^p(\mu_k, \mu) > n^{\alpha} x
\right )  \ll  \frac{n^{1-\alpha}}{x} \int_0^\infty t^{p-1} H(t) {\bf 1}_{t >n^{\alpha/p}} dt  + 
\frac{n^{1-q \alpha}}{x^q} 
 \left( \int_0^\infty t^{p-1} (H(t))^{1/q} {\bf 1}_{t \leq n^{\alpha/p}} dt
 \right )^q \, .
$$
Hence, it remains to prove that 
\begin{equation}\label{2series}
\sum_{n=1}^{\infty} 
n^{\alpha( r -1) -1}  \int_0^\infty t^{p-1} H(t) {\bf 1}_{t^{p/\alpha} >n} dt < \infty \quad \text{and} \quad 
\sum_{n=1}^{\infty} 
n^{\alpha( r -q) -1}  \left( \int_0^\infty t^{p-1} (H(t))^{1/q} {\bf 1}_{t^{p/\alpha} \leq n} dt
 \right )^q < \infty \, .
\end{equation}
Interverting the sum and the integral, we easily get that 
\begin{equation*}
\sum_{n=1}^{\infty} 
n^{\alpha( r -1) -1}  \int_0^\infty t^{p-1} H(t) {\bf 1}_{t^{p/\alpha} >n} dt \ll \int_0^\infty t^{pr-1} H(t) dt  \ll \|X \|_{pr}^{pr} < \infty \, .
\end{equation*}
Arguing as in \eqref{Heq} with $\beta <(q-1)/q$, we get 
$$
  \left( \int_0^\infty t^{p-1} (H(t))^{1/q} {\bf 1}_{t^{p/\alpha} \leq n} dt   \right )^q  \ll
n^{\alpha(q-1-q\beta)/p}  \int_0^{\infty} t^{q(p-1 + \beta ) } H (t){\bf 1}_{t^{p/\alpha} \leq n} dt   \, .
$$
Hence, the second series in \eqref{2series} will be summable provided
\begin{equation}\label{ser2}
\sum_{n=1}^{\infty}
n^{\alpha( r -q) + \alpha(q-1-q\beta)/p -1}   \int_0^{\infty} t^{q(p-1 + \beta ) } H (t){\bf 1}_{t^{p/\alpha} \leq n} dt  < \infty
\end{equation}
Taking $\beta$ close enough to $(q-1)/q$ so that $\alpha( r -q) + \alpha(q-1-q\beta)/p<0$ and
interverting the sum and the integral, we  get that 
\begin{equation*}
\sum_{n=1}^{\infty} 
n^{\alpha( r -q) + \alpha(q-1-q\beta)/p -1}   \int_0^{\infty} t^{q(p-1 + \beta ) } H (t){\bf 1}_{t^{p/\alpha} \leq n} dt 
 \ll \int_0^\infty t^{pr-1} H(t) dt  \ll \|X \|_{pr}^{pr} < \infty \, ,
 \end{equation*}
which ends  the proof of \eqref{2series} and then the proof of the theorem when $p > d(r-1)/r$.


\smallskip

\noindent $\bullet$  If $p<d(r-1)/r$, let $q \in (r,2]$, $s \in (1,r)$ such that $p<d(s-1)/s$, and let  $M=n^{(p-d(1-\alpha))/(dp)}$. From \eqref{start}, 
\eqref{Abound3} and \eqref{Bbound3}, we get the upper bound
\begin{multline}\label{bkbis}
{\mathbb P}\left ( \max_{1 \leq k \leq n} k W_p^p(\mu_k, \mu) > n^{\alpha}  x
\right )  \ll  \frac{n^{s-s \alpha}}{n^{sp/d} x^s} 
 \left( \int_0^\infty t^{p-1} (H(t))^{1/s} {\bf 1}_{t > n^{(p-d(1-\alpha))/(dp)}} dt
 \right )^s \\
+ 
\frac{n^{q-q \alpha}}{n^{qp/d} x^q} 
 \left( \int_0^\infty t^{p-1} (H(t))^{1/q} {\bf 1}_{t \leq n^{(p-d(1-\alpha))/(dp)}} dt
 \right )^q \, .
\end{multline}
Proceeding as in \eqref{ser2} (taking the quantity $\beta < (q-1)/q$ close enough to $(q-1)/q$ in such a way that 
$
(p-d(1- \alpha)) ((r-q) + (q-1- \beta q)/p )<0
$), 
we get that 
\begin{equation}\label{bkbis2}
\sum_{n=1}^\infty n^{(pr-(1-\alpha)rd-d)/d} \frac{n^{q-q \alpha}}{n^{qp/d} }  \left( \int_0^\infty t^{p-1} (H(t))^{1/q} {\bf 1}_{t \leq n^{(p-d(1-\alpha))/(dp)}} dt
 \right )^q \ll \|X \|_{pr}^{pr} < \infty \, .
\end{equation}
In the same, we get 
\begin{equation}\label{bkbis3}
\sum_{n=1}^\infty n^{(pr-(1-\alpha)rd-d)/d} \frac{n^{s-s \alpha}}{n^{sp/d} } 
 \left( \int_0^\infty t^{p-1} (H(t))^{1/s} {\bf 1}_{t > n^{(p-d(1-\alpha))/(dp)}} dt
 \right )^s  \ll \|X \|_{pr}^{pr} < \infty \, .
\end{equation}
The second item of Theorem \ref{BK1} follows from \eqref{bkbis}, \eqref{bkbis2}
and \eqref{bkbis3}.

\smallskip

\noindent $\bullet$  If $p<d(r-1)/r$, let $q \in (r,2]$, $s \in (1,r)$ such that $p<d(s-1)/s$, and let  $M=(\log n )^{1/pr}$. From \eqref{start}, 
\eqref{Abound3} and \eqref{Bbound3}, we get the upper bound
\begin{multline}\label{bkter}
{\mathbb P}\left ( \max_{1 \leq k \leq n} k W_p^p(\mu_k, \mu) > n^{(d-p)/d} (\log n)^{1/r} x
\right )  \ll  \frac{1}{(\log n)^{s/r} x^s} 
 \left( \int_0^\infty t^{p-1} (H(t))^{1/s} {\bf 1}_{t > (\log n )^{1/pr}} dt
 \right )^s \\
+ 
 \frac{1}{(\log n)^{q/r} x^q}
 \left( \int_0^\infty t^{p-1} (H(t))^{1/q} {\bf 1}_{t \leq  (\log n )^{1/pr}} dt
 \right )^q \, .
\end{multline}
Proceeding as in \eqref{ser2} (taking the quantity $\beta < (q-1)/q$ close enough to $(q-1)/q$ in such a way that  
$
(q/r)- (q-1-\beta q)/(pr) >1
$), 
we get that 
\begin{equation}\label{bkter2}
\sum_{n=1}^\infty \frac{1}{n(\log n)^{q/r}}  \left( \int_0^\infty t^{p-1} (H(t))^{1/q} {\bf 1}_{t \leq (\log n )^{1/pr}} dt
 \right )^q \ll \|X \|_{pr}^{pr} < \infty \, .
\end{equation}
In the same, we get 
\begin{equation}\label{bkter3}
\sum_{n=1}^\infty \frac{1}{n(\log n)^{s/r}} 
 \left( \int_0^\infty t^{p-1} (H(t))^{1/s} {\bf 1}_{t >  (\log n )^{1/pr}} dt
 \right )^s  \ll \|X \|_{pr}^{pr} < \infty \, .
\end{equation}
The third item of Theorem \ref{BK1} follows from \eqref{bkter}, \eqref{bkter2}
and \eqref{bkter3}.

\subsection{Proof of Theorem \ref{BK2}}
Let $r >2$. As in the proof of Theorem \ref{F-Nag}, we assume without loss of generality that $\|X\|_{rp, w}=1$; hence, we can use directly some of  the upper bounds given in the proof 
of Theorem \ref{F-Nag}. 

From \eqref{fgu1}, we see that
\begin{equation}\label{fgu1bis}
\max_{1 \leq k \leq n} k W_p^p(\mu_k, \mu) \leq C \sum_{m \geq 0}  2^{pm}   \max_{1 \leq k \leq n} \left | k\mu_k(B_m)- k\mu(B_m) \right | + C  \max_{1 \leq k \leq n} kV_k^p \, .
\end{equation}
Now, for any $x>0$
$$
{\mathbb P}\left ( \max_{1 \leq k \leq n} kV_k^p > x /(2C) \right )  \leq \sum_{k=1}^n {\mathbb P}\left (  kV_k^p > x/(2C)  \right ) \leq n \max_{1 \leq k \leq n}{\mathbb P}\left ( 
 kV_k^p > x /(2C)  \right )\, .
$$
By \eqref{fgu2}, it follows that, for any $x>0$, 
\begin{equation}\label{FGU3}
{\mathbb P}\left ( \max_{1 \leq k \leq n} kV_k^p > x /(2C) \right ) \leq n \max_{1 \leq k \leq n} a(k,x/k) \leq na(n,x/n) \, ,
\end{equation}
the last inequality being true because $k \rightarrow a(k,x/k)$ is increasing. Now, by definition of $a(n,x)$, we infer that, for any 
$$
 \alpha \in \left ( \max \left ( \frac 1 2, \frac {d-p}{d} \right ) ,  1
 \right ] \, ,
$$
\begin{equation*}
\sum_{n=1}^{\infty} 
n^{\alpha r -2} {\mathbb P}\left ( \max_{1 \leq k \leq n} k V_k^p > n^{\alpha} x/(2C)
\right )  < \infty  \, .
\end{equation*}
Hence, it remains to prove that 
\beq \label{Th34-p1-15-dec}
\sum_{n=1}^{\infty} 
n^{\alpha r -2}
{\mathbb P} \left ( \sum_{m \geq 0}  2^{pm}  \max_{1 \leq k \leq n} \left |k \mu_k(B_m)- k\mu(B_m) \right | \geq  n^{\alpha} x/(2C) \right )  < \infty \, .
\eeq
Arguing as in the proof of Theorem \ref{F-Nag}, and using a maximal version of Rosenthal's inequality (see for instance \cite{P}), we get that, for any  $q>r$ and $M>0$, 
\begin{multline}\label{3series}
{\mathbb P} \left ( \sum_{m \geq 0}  2^{pm}  \max_{1 \leq k \leq n} \left |k \mu_k(B_m)- k\mu(B_m) \right | \geq  n^{\alpha} x/(2C) \right ) 
\ll  \frac{n^{(1-\alpha)q}}{x^{q}n^{q/2}}  \left ( \int_0^{\infty} t^{p-1} \sqrt{H (t)}  dt  \right )^q \\
+  \frac{n^{1-q\alpha}}{x^q}  \left (\int_0^{\infty} t^{p-1} H^{1/q} (t) {\bf 1}_{t \leq M} dt \right)^q +  \frac{n^{1-2\alpha}}{ x^2} \left ( \int_0^{\infty} t^{p-1} \sqrt{H (t)} {\bf 1}_{t > M} dt \right )^2 \, .
\end{multline}
Clearly, since $ \alpha \in (1/2,1]$, taking $q$ large enough, we get that 
\begin{equation*}
\sum_{n=1}^{\infty} 
n^{\alpha r -2} \frac{n^{(1-\alpha)q}}{n^{q/2}}  \left ( \int_0^{\infty} t^{p-1} \sqrt{H (t)}  dt  \right )^q < \infty \, .
\end{equation*}
Let $M=n^{\alpha/p}$ and  $\beta >1/2$. Applying H\"older's inequality, we get 
$$
  \left( \int_0^\infty t^{p-1} \sqrt{H(t)} {\bf 1}_{t^{p/\alpha} > n} dt   \right )^2  \ll
n^{\alpha(1-2 \beta)/p}  \int_0^{\infty} t^{2(p-1 + \beta ) } H (t){\bf 1}_{t^{p/\alpha} > n} dt   \, .
$$
Hence, the sum over $n$ of the last term in \eqref{3series} multiplied by $n^{\alpha r -2}$ will be finite provided
\begin{equation*}
\sum_{n=1}^{\infty}
n^{\alpha( r -2) + \alpha(1-2\beta)/p -1}   \int_0^{\infty} t^{2(p-1 + \beta ) } H (t){\bf 1}_{t^{p/\alpha} > n} dt  < \infty \, .
\end{equation*}
Taking $\beta$ close enough to $1/2$ so that $\alpha( r -2) + \alpha(1-2\beta)/p>0$ and
interverting the sum and the integral, we  get that 
\begin{equation*}
\sum_{n=1}^{\infty} 
n^{\alpha( r -2) + \alpha(1-2\beta)/p -1}   \int_0^{\infty} t^{2(p-1 + \beta ) } H (t){\bf 1}_{t^{p/\alpha} > n} dt 
 \ll \int_0^\infty t^{pr-1} H(t) dt  \ll \|X \|_{pr}^{pr} < \infty \, .
 \end{equation*}
 
Arguing as in \eqref{Heq} with $\beta <(q-1)/q$, we get 
$$
  \left( \int_0^\infty t^{p-1} H^{1/q}(t) {\bf 1}_{t^{p/\alpha} \leq  n} dt   \right )^q  \ll
n^{\alpha(q-1- \beta q)/p}  \int_0^{\infty} t^{q(p-1 + \beta ) } H (t){\bf 1}_{t^{p/\alpha} \leq  n} dt   \, .
$$
Hence,  the sum over $n$ of the second term in \eqref{3series} multiplied by $n^{\alpha r -2}$ will be finite provided
\begin{equation*}
\sum_{n=1}^{\infty}
n^{\alpha( r -q) + \alpha(q-1-\beta q)/p -1}   \int_0^{\infty} t^{q(p-1 + \beta ) } H (t){\bf 1}_{t^{p/\alpha} \leq  n} dt  < \infty \, .
\end{equation*}
Taking $\beta$ close enough to $(q-1)/q$ so that $\alpha( r -q) + \alpha(q-1-\beta q)/p<0$ and
interverting the sum and the integral, we  get that 
\begin{equation*}
\sum_{n=1}^{\infty} 
n^{\alpha( r -q) + \alpha(q-1- \beta q)/p -1}   \int_0^{\infty} t^{q(p-1 + \beta ) } H (t){\bf 1}_{t^{p/\alpha} \leq  n} dt 
 \ll \int_0^\infty t^{pr-1} H(t) dt  \ll \|X \|_{pr}^{pr} < \infty \, .
 \end{equation*}
 All the previous considerations end the proof of  \eqref{Th34-p1-15-dec} and then of the theorem.

\subsection{Proof of Theorem \ref{LILBornee}}
Recall that 
\[
{\mathcal D}_p ({\mathcal R}_{B_m}\mu_n,{\mathcal R}_{B_m} \mu)  = \frac{2^p -1}{2} \sum_{\ell \geq 1} 2^{-p \ell} \sum_{F \in {\mathcal P}_{\ell}}
 \left |  \frac{\mu_n ({\tilde F}_m)}{\mu_n (B_m) } -  \frac{\mu ({\tilde F}_m)}{\mu (B_m) }\right | \, , 
\]
where ${\tilde F}_m = 2^m F \cap B_m$ (see \eqref{defclaDpR}). Define, for any $k \geq 1$, 
\[
n_k = [{\rm e}^{k^{1/2}}] \text{ and } m_k = n_{k+1} -n_k \, .
\]
Note that $\ n_{k+1} \leq 2 {\rm e} \, n_k$ and $m_k \sim 2^{-1} k^{- 1/2}n_k$, as $k \rightarrow \infty$. Setting 
\[
\mu_{n_k,n} = \frac{1}{n-n_k} \sum_{i=n_k+1}^n  \delta_{X_i} \, , 
\]
we first write that, for $n_k +1 \leq n < n_{k+1}$, 
\begin{multline*}
 \frac{\mu_n ({\tilde F}_m)}{\mu_n (B_m) } -  \frac{\mu ({\tilde F}_m)}{\mu (B_m) }   \\
  =  \frac{ n_k ( \mu_{n_k}({\tilde F}_m)- \mu ({\tilde F}_m) ) }{n \mu_n (B_m) } +  \frac{(n-n_k)  (\mu_{n_k,n} ({\tilde F}_m)- \mu ({\tilde F}_m) ) }{n \mu_n (B_m) } +  \frac{\mu(B_m ) - \mu_n(B_m) }{\mu_n (B_m) \mu (B_m) }  \mu ({\tilde F}_m)  \, .
\end{multline*}
Taking into account that, for any positive measure $\nu$,
\[
 \sum_{\ell \geq 1} 2^{-p \ell} \sum_{F \in {\mathcal P}_{\ell}} \nu ({\tilde F}_m) \leq (2^p-1)^{-1}\nu ( B_m)  \, , 
\]
simple algebras lead to the following inequality:  for $n_k +1 \leq n < n_{k+1}$, 
\begin{multline*}
 \sum_{\ell \geq 1} 2^{-p \ell} \sum_{F \in {\mathcal P}_{\ell}}   \frac{ n_k | \mu_{n_k}({\tilde F}_m)- \mu ({\tilde F}_m) | }{n \mu_n (B_m) } \leq    \sum_{\ell \geq 1} 2^{-p \ell} \sum_{F \in {\mathcal P}_{\ell}}   \left | \frac{  \mu_{n_k}({\tilde F}_m) }{ \mu_{n_k} (B_m) }  - \frac{ \mu({\tilde F}_m) }{ \mu(B_m) }   \right |  \\ + \frac{1}{2^p-1}  \frac{| \mu_{n_k} (B_m) - \mu_n (B_m) | }{ \mu_n (B_m) } 
+ \frac{1}{2^p-1}  \frac{| \mu_{n_k} (B_m) - \mu (B_m) | }{ \mu_n (B_m) }   \, .
\end{multline*}
Similarly,  for $n_k +1 \leq n < n_{k+1}$, 
\begin{multline*}
 \sum_{\ell \geq 1} 2^{-p \ell} \sum_{F \in {\mathcal P}_{\ell}}   \frac{ (n-n_k ) | \mu_{n_k,n}({\tilde F}_m)- \mu ({\tilde F}_m) | }{n \mu_n (B_m) } \leq \frac{(n-n_k)}{n}  \sum_{\ell \geq 1} 2^{-p \ell} \sum_{F \in {\mathcal P}_{\ell}}   \left | \frac{  \mu_{n_k,n}({\tilde F}_m) }{ \mu_{n_k,n} (B_m) }  - \frac{ \mu({\tilde F}_m) }{ \mu(B_m) }   \right |  \\ + \frac{ (n-n_k)  }{(2^p-1) n }  \frac{| \mu_{n_k,n} (B_m) - \mu_n (B_m) | }{ \mu_n (B_m) } 
+ \frac{ (n-n_k)}{(2^p-1) n }  \frac{| \mu_{n_k,n} (B_m) - \mu (B_m) | }{ \mu_n (B_m) }   \, .
\end{multline*}
So overall,  for $n_k +1 \leq n < n_{k+1}$, 
\begin{multline} \label{ine1-13-07}
{\mathcal D}_p ({\mathcal R}_{B_m}\mu_n,{\mathcal R}_{B_m} \mu)  \leq {\mathcal D}_p ({\mathcal R}_{B_m}\mu_{n_k},{\mathcal R}_{B_m} \mu)  + \frac{(n-n_k)}{n}   {\mathcal D}_p ({\mathcal R}_{B_m}\mu_{n_k,n},{\mathcal R}_{B_m} \mu) \\
+  \frac{1}{2}  \frac{| \mu_{n_k} (B_m) - \mu_n (B_m) | }{ \mu_n (B_m) } 
+ \frac{1}{2}  \frac{| \mu_{n_k} (B_m) - \mu (B_m) | }{ \mu_n (B_m) } +  \frac{1}{2}  \frac{| \mu_{n} (B_m) - \mu (B_m) | }{ \mu_n (B_m) }   \\
+ \frac{ (n-n_k) }{2 n } \frac{| \mu_{n_k,n} (B_m) - \mu_n (B_m) | }{ \mu_n (B_m) } 
+ \frac{ (n-n_k) }{ 2 n }   \frac{| \mu_{n_k,n} (B_m) - \mu (B_m) | }{ \mu_n (B_m) } \, .
\end{multline}

\smallskip

\noindent $\bullet$ If $p >d/2$, 
let 
\[
v_n =  \sqrt{\frac{\log \log n}{ n } } \quad  \text{and}  \quad  V= \int_0^{\infty} t^{p-1} \sqrt{H(t)} dt \, .
\]
Starting from \eqref{defclaDp} and  considering \eqref{ine1-13-07}, it follows that 
\begin{multline} \label{decpreuveLIL}
 \max_{n_k+1 \leq n \leq n_{k+1}}\frac{{\mathcal D}_p (\mu_n, \mu) }{v_n} 
 \leq 
 \sum_{m \geq 0} 2^{pm} \mu(B_m)  \frac{ {\mathcal D}_p ({\mathcal R}_{B_m}\mu_{n_k},{\mathcal R}_{B_m} \mu)  }{v_{n_{k+1}}}  \\ +  \max_{n_k+1 \leq n \leq n_{k+1}}  \frac{(n-n_k)}{n}   \sum_{m \geq 0} 2^{pm} \mu(B_m)  \frac{ {\mathcal D}_p ({\mathcal R}_{B_m}\mu_{n_k,n},{\mathcal R}_{B_m} \mu)  }{v_{n_{k+1}}}\\
+  \max_{n_k \leq n \leq n_{k+1}}   \Big ( 2+\frac{( n-n_k) }{ 2n }   \Big )   \sum_{m \geq 0} 2^{pm}   \frac{| \mu_{n} (B_m) - \mu (B_m) | }{v_{n_{k+1}} } \\
+   \max_{n_k +1 \leq  n \leq n_{k+1}}   \frac{ ( n-n_k) }{  n }    \sum_{m \geq 0} 2^{pm}   \frac{| \mu_{n_k, n} (B_m) - \mu (B_m) | }{v_{n_{k+1}} }  \, .
\end{multline}
We first deal with the third term in the right-hand side of \eqref{decpreuveLIL}. With this aim, let 
\beq \label{choiceM-14-dec}
M=M_{n_k}  = \delta \left ( \frac{n_k}{LL n_k }\right )^{1/(2p)}
\eeq
with $\delta$ a positive constant not depending on $n_k$ that will be chosen later, and the notations $Lx = \log (x \vee e)$ and $ LLx = L(Lx)$. Define ${\mathcal C}_{n_k} = [-M_{n_k} , M_{n_k} ]^d$ and note that 
\[
| \mu_{n} (B_m) - \mu (B_m) |  \leq  | \mu_{n} (B_m \cap {\mathcal C}_{n_k}) - \mu (B_m\cap {\mathcal C}_{n_k}) | + | \mu_{n} (B_m \cap {\mathcal C}^c_{n_k}) - \mu (B_m\cap {\mathcal C}^c_{n_k}) |
\]
Clearly
\beq \label{p1-14-dec-18}
  \sum_{m \geq 0} 2^{pm}  \frac{\mu (B_m\cap {\mathcal C}^c_{n_k})}{v_{n_{k+1}} }  \ll  \frac{1}{v_{n_{k+1}} } \int_0^{\infty}  t^{p-1} H(t)  {\bf 1}_{t > M_{n_k}  }dt \ll 
  \int_0^{\infty}  t^{2p-1} H(t)  {\bf 1}_{t > M_{n_k}  }dt \rightarrow 0 \, , \, \text{ as $k \rightarrow \infty$ \, .}
\eeq
On the other hand, 
\[
 \max_{n_k \leq n \leq n_{k+1}}   \mu_{n} (B_m \cap {\mathcal C}^c_{n_k})  \leq  \frac{1}{n_k} \sum_{i=1}^{n_{k+1}} {\bf 1}_{\{ X_i \in B_m\}} {\bf 1}_{\{ |X_i|  >M_{n_k} \}} 
 \ll  \frac{1}{n_k} \sum_{i=1}^{n_{k+1}} {\bf 1}_{\{ X_i \in B_m\}} {\bf 1}_{\{ |X_i|  > c M_i\}}  \, ,
\]
where $c$ is a universal positive constant.  Hence to prove that 
\beq \label{p2-14-dec-18}
 \max_{n_k \leq n \leq n_{k+1}}   \sum_{m \geq 0} 2^{pm}   \frac{ \mu_{n} (B_m \cap {\mathcal C}^c_{n_k})  }{v_{n_{k+1}} }   \rightarrow 0 \, , \, \text{ almost surely, as $k \rightarrow \infty$,}
\eeq
it suffices to prove that 
\beq \label{p2-15-dec-18bis}
\frac{1}{ n v_{n} } \sum_{i=1}^{n}  \sum_{m \geq 0} 2^{pm}  {\bf 1}_{\{ X_i \in B_m\}} {\bf 1}_{\{ |X_i|  > c M_i\}}   \rightarrow 0 \, , \, \text{ almost surely, as $n \rightarrow \infty$.}
\eeq
But
\begin{multline*}
\sum_{i \geq 1}  \frac{1}{ \sqrt{ i  LL i }}\sum_{m \geq 0} 2^{pm}    \P ( X_i \in B_m,  |X_i|  > c M_{i} ) \ll   \sum_{i \geq 1}  \frac{1}{ \sqrt{ i  LL i }} \int_0^{\infty}  t^{p-1} H(t)  {\bf 1}_{t > c M_i }dt  \\ \ll \int_0^{\infty}  t^{2p-1} H(t) dt < \infty \, .
 \end{multline*}
By using Kronecker's lemma and recalling that $n v_n = (n LL n)^{1/2}$,  this shows that \eqref{p2-15-dec-18bis} holds and so \eqref{p2-14-dec-18} does also. 

We show now that there exists a positive constant $C$ such that, almost surely, 
\beq \label{p3-14-dec-18}
\limsup_{k \rightarrow \infty}    \max_{n_k \leq n \leq n_{k+1}}   \sum_{m \geq 0} 2^{pm}   \frac{\big |  \mu_{n} (B_m \cap {\mathcal C}_{n_k})  -   \mu (B_m \cap {\mathcal C}_{n_k}) \big |  }{v_{n_{k+1}} }  \leq CV    \, .
\eeq
Using Markov's inequality and next  Rosenthal's inequality   (with the constants given in (4.2) of Theorem 4.1  in \cite{P}), as in the proof of Theorem \ref{F-Nag}, we infer that there exist positive universal constants $c_1$ and  $c_2$  such that for any $q > 2$ and  $\lambda >0$,  
\begin{multline*}
\P \left (   \max_{n_k \leq n \leq n_{k+1}}   \sum_{m \geq 0} 2^{pm}  | \mu_{n} (B_m\cap {\mathcal C}_{n_k}) - \mu (B_m\cap {\mathcal C}_{n_k}) |   
\geq  \lambda  V v_{n_{k+1}}  \right ) \\  \leq
  \left ( \frac{ c_1 }{  \lambda V  n_k v_{n_{k+1}} }  \right )^q q^{q/2}  n_{k+1}^{q/2} V^q +  \left ( \frac{ c_2 }{  \lambda V  n_k v_{n_{k+1}} }  \right )^q q^{q}  n_{k+1}  
  \left (\int_0^{\infty}  t^{p-1} (H(t) )^{1/q} {\bf 1}_{t \leq M_{n_k} }dt  \right )^q   \, .
 \end{multline*}
 Select now 
 \[
 q=q_k=\gamma \log \log n_k \  \text{  with $\gamma >2$}, 
 \]
 and take $\lambda= \lambda_{\gamma} = 2 c_1 {\rm e}^2 \sqrt{\gamma} $.  With this choice of $q_k$, it follows that 
\[
\sum_{k \geq 1}   \left ( \frac{ c_1 }{  \lambda_{\gamma}  n_k v_{n_{k+1}} }  \right )^{q_k} q_k^{q_k/2}  n_{k+1}^{q_k/2} \leq  \sum_{k \geq 1}  
 \left ( \frac{  2 c_1 {\rm e} \sqrt{\gamma}    }{  \lambda_{\gamma}  }   \right )^{q_k}   = \sum_{k \geq 1}  
  {\rm e}^{-q_k}   < \infty \, .
\]
On the other hand, by H\"older's inequality, setting $\beta = p-2p/q$, 
\[
 \left (\int_0^{\infty}  t^{p-1} (H(t) )^{1/q} {\bf 1}_{t \leq M_{n_k}  } dt  \right )^q \leq p^{-q}2^{q-1}  + 2^{q-1} M_{n_k}^{pq -2p} \beta^{1-q}  \int_0^{\infty}  t^{2p-1}H(t) dt  \, .
\]
Concerning the constant $\delta$ appearing in the selection of $M_{n_k}$ given in \eqref{choiceM-14-dec}, select it such that 
\[
\frac{4 \sqrt{2 {\rm e}}c_2 \delta^p \gamma }{  \lambda_{\gamma} V p} = {\rm e}^{-1} \, .
\]
Let  $K_1 $ be such that $q_{K_1} \geq 4$.  It follows that 
\begin{multline*}
\sum_{k \geq K_1}   \left ( \frac{ c_2 }{  \varepsilon V  n_k v_{n_{k+1}} }  \right )^q q^{q}  n_{k+1}  \left (\int_0^{\infty}  t^{p-1} (H(t) )^{1/q} {\bf 1}_{t \leq M_{n_k}  }dt  \right )^q  \\
\ll  \delta^{-2p}\sum_{k \geq K_1}   \left ( \frac{4 \sqrt{2 {\rm e}}c_2 \delta^p \gamma }{  \lambda_{\gamma} V p}   \right )^{q_k}  n_{k+1}  \frac{\log \log n_k} {n_k} \ll
\delta^{-2p} \sum_{k \geq K_1} 
{\rm  e}^{-q_k}  \log \log n_k < \infty  \, .
 \end{multline*}
So, overall, 
\[
\sum_{k \geq K_1} \P \left (   \max_{n_k \leq n \leq n_{k+1}}   \sum_{m \geq 0} 2^{pm}  | \mu_{n} (B_m\cap {\mathcal C}_{n_k}) - \mu (B_m\cap {\mathcal C}_{n_k}) |    \geq  \lambda_{\gamma}  V v_{n_{k+1}}  \right ) < \infty  \, ,
\]
which proves \eqref{p3-14-dec-18} with $  C=\lambda_{\gamma} $ by the direct part of the  Borel-Cantelli lemma.  Hence combining \eqref{p1-14-dec-18}, \eqref{p2-14-dec-18} and \eqref{p3-14-dec-18}, it follows that, almost surely, 
\beq \label{p4-14-dec-18}
\limsup_{k \rightarrow \infty}  \max_{n_k \leq n \leq n_{k+1}}    \sum_{m \geq 0} 2^{pm}   \frac{\big |  \mu_{n} (B_m )  -   \mu (B_m ) \big |  }{v_{n_{k+1}} }  \leq \lambda_{\gamma} V    \, .
\eeq
With similar arguments, one can prove that, almost surely, 
\beq \label{p5-14-dec-18}
\limsup_{k \rightarrow \infty}  \max_{n_k \leq n \leq n_{k+1}}  \frac{n-n_k}{n}  \sum_{m \geq 0} 2^{pm}   \frac{\big |  \mu_{n_k,n} (B_m )  -   \mu (B_m ) \big |  }{v_{n_{k+1}} }  =0     \, .
\eeq
It follows that, almost surely, 
\begin{multline} \label{ine2-13-07}
\limsup_{k \rightarrow \infty}  \max_{n_k+1 \leq n \leq n_{k+1}}\frac{{\mathcal D}_p (\mu_n, \mu) }{v_n} 
 \leq 2 \lambda_{\gamma} V + \limsup_{k \rightarrow \infty}  
 \sum_{m \geq 0} 2^{pm} \mu(B_m)  \frac{ {\mathcal D}_p ({\mathcal R}_{B_m}\mu_{n_k},{\mathcal R}_{B_m} \mu)  }{v_{n_{k+1}}}  \\ + \limsup_{k \rightarrow \infty}    \max_{n_k+1 \leq n \leq n_{k+1}}   \frac{(n-n_k)}{n}  \sum_{m \geq 0} 2^{pm} \mu(B_m)  \frac{ {\mathcal D}_p ({\mathcal R}_{B_m}\mu_{n_k,n},{\mathcal R}_{B_m} \mu)  }{v_{n_{k+1}}}  \, .
\end{multline}
Let now
\[
s_k = \left  [ \frac{  k^{1/2}}{p\ln 2} \right ] \, .
\]
Note that 
\begin{multline} \label{ine3ante-13-07}
 \sum_{m \geq s_k+2} 2^{pm} \mu(B_m)  {\mathcal D}_p ({\mathcal R}_{B_m}\mu_{n_k},{\mathcal R}_{B_m} \mu)
\leq  \sum_{m \geq s_k+2} 2^{pm} \mu ( B_m )  \\ \leq {\tilde C}_p  \int_{ 2^{s_k}}^{\infty}  t^{p-1} H(t) dt \leq  {\tilde C}_p  2^{-p s_k} \int_{ 0}^{\infty}  t^{2p-1} H(t) dt  \, .
\end{multline}
It follows that 
\beq \label{ine3-13-07}
\lim_{k \rightarrow \infty} \sum_{m \geq s_k+1} 2^{pm} \mu(B_m) \frac{ {\mathcal D}_p ({\mathcal R}_{B_m}\mu_{n_k},{\mathcal R}_{B_m} \mu)  }{v_{n_{k+1}} }  = 0 \quad a.s.  
\eeq
Next, let 
\[
b_m =  \int_{ 2^{m-2}}^{2^{m-1}}  t^{p-1} \sqrt{H(t)} dt {\bf 1}_{m \geq 2} +   \int_{ 0}^{1}  t^{p-1} \sqrt{H(t)} dt {\bf 1}_{m \leq 2}  \  \text{ and }  \ B= \sum_{m \geq 0} b_m = V +  \int_{ 0}^{1}  t^{p-1} \sqrt{H(t)} dt \, .
\]
Note that 
\begin{multline*}
\P \left (  \sum_{m =0}^{ s_k+1} 2^{pm} \mu(B_m)  {\mathcal D}_p ({\mathcal R}_{B_m}\mu_{n_k},{\mathcal R}_{B_m} \mu)  \geq  C B  v_{n_{k+1}}  \right )  \\
\leq   \sum_{m =0}^{ s_k+1}  \P \left ( 2^{pm} \mu(B_m)  {\mathcal D}_p ({\mathcal R}_{B_m}\mu_{n_k},{\mathcal R}_{B_m} \mu)  \geq    C b_m  v_{n_{k+1}}  \right ) \, .
\end{multline*}
Proceeding as in the proof of Theorem 2 in \cite{FG} (case $p>d/2$), and noting that 
\begin{multline} \label{evident}
\mu(B_m)  \leq \P ( |X| > 2^{m-1} ) \leq \left (  \frac{1}{2^{m-2}} \int_{ 2^{m-2}}^{2^{m-1}}  \sqrt{H(t)} dt \right )^2 \\
 \leq \left (  \frac{1}{2^{(m-2)p}} \int_{ 2^{m-2}}^{2^{m-1}} t^{p-1} \sqrt{H(t)} dt \right )^2 = 2^{4p} 2^{-2mp} b_m^2  \, ,
\end{multline}
we derive that, there exists a positive universal constant $a$ such that 
\begin{multline*}
\P \left (  \sum_{m =0}^{ s_k+1} 2^{pm} \mu(B_m)  {\mathcal D}_p ({\mathcal R}_{B_m}\mu_{n_k},{\mathcal R}_{B_m} \mu)  \geq  C B  v_{n_{k+1}}  \right )  \\
\leq   \sum_{m =0}^{ s_k+1}  \exp \left (   - a C^2 b^2_m n_k  v^2_{n_k} / (2^{2pm} \mu ( B_m))   \right ) \\
 \leq  (s_k+2)  \exp \left (   - \frac{a C^2 }{2^{4p}}  \log \log n_k   \right ) \ll  \frac{k^{1/2}}{ k^{aC^2 /2^{4p+1}}}  \, .
\end{multline*}
Therefore, if $C$ is large enough, 
\beq \label{ine4-13-07bis}
\sum_{k \geq 1} \P \left (  \sum_{m =0}^{ s_k+1} 2^{pm} \mu(B_m)  {\mathcal D}_p ({\mathcal R}_{B_m}\mu_{n_k},{\mathcal R}_{B_m} \mu)  \geq  C B  v_{n_{k+1}}  \right )  < \infty\, .
\eeq
Starting from \eqref{ine3-13-07} and  \eqref{ine4-13-07bis}, it follows that 
\beq \label{ine5-13-07}
\limsup_{k \rightarrow \infty} \sum_{m  \geq 0} 2^{pm} \mu(B_m) \frac{ {\mathcal D}_p ({\mathcal R}_{B_m}\mu_{n_k},{\mathcal R}_{B_m} \mu)  }{v_{n_{k+1}} }  \leq CB  \quad a.s.  
\eeq
On another hand, using \eqref{ine3ante-13-07}, we get that 
\beq \label{ine6-13-07}
\limsup_{k \rightarrow \infty}  \max_{n_k+1 \leq n \leq n_{k+1}}   \frac{(n-n_k)}{n}  \sum_{m \geq s_k +2} 2^{pm} \mu(B_m)  \frac{ {\mathcal D}_p ({\mathcal R}_{B_m}\mu_{n_k,n},{\mathcal R}_{B_m} \mu)  }{v_{n_{k+1}}} = 0 \quad a.s.  
\eeq
In addition, noticing that $ {\mathcal D}_p ({\mathcal R}_{B_m}\mu_{n_k,n},{\mathcal R}_{B_m} \mu) =^{\mathcal D}  {\mathcal D}_p ({\mathcal R}_{B_m}\mu_{n-n_k},{\mathcal R}_{B_m} \mu) $, we get that
\begin{multline*}
\P \left (  \max_{n_k+1 \leq n \leq n_{k+1}}   \frac{(n-n_k)}{n}  \sum_{m = 0}^{s_k+1} 2^{pm} \mu(B_m)  {\mathcal D}_p ({\mathcal R}_{B_m}\mu_{n_k,n},{\mathcal R}_{B_m} \mu)   \geq  C B  v_{n_{k+1}}  \right )  \\
\leq  \sum_{n=n_k+1}^{ n_{k+1}}  \sum_{m =0}^{ s_k+1}   \P \left (  \frac{(n-n_k)}{n}  2^{pm} \mu(B_m)  {\mathcal D}_p ({\mathcal R}_{B_m}\mu_{n-n_k},{\mathcal R}_{B_m} \mu)  \geq    C  b_m  v_{n_{k+1}}  \right )  \, .  \\
\end{multline*}
Proceeding as before,   we get  that
\begin{multline} \label{ine4-13-07-bisbis}
\sum_{k \geq 4}\P \left (  \max_{n_k+1 \leq n \leq n_{k+1}}   \frac{(n-n_k)}{n}  \sum_{m = 0}^{s_k+1} 2^{pm} \mu(B_m)  {\mathcal D}_p ({\mathcal R}_{B_m}\mu_{n_k,n},{\mathcal R}_{B_m} \mu)   \geq  C  B  v_{n_{k+1}}  \right )  \\
\ll \sum_{k \geq 4}   \sum_{n=n_k+1}^{ n_{k+1}}    \sum_{m =0}^{ s_k+1}  \exp \left (   - \frac{a C^2 n^2 b^2_m   v^2_{n_k} }{(n-n_k) 2^{2pm} \mu ( B_m)}   \right )  \\ 
\ll \sum_{k \geq 4}   ( n_{k+1} -n_k)    \sum_{m =0}^{ s_k+1}  \exp \left (   - \frac{a C^2 n_k^2 b^2_m   v^2_{n_k} }{(n_{k+1}-n_k) 2^{2pm} \mu ( B_m)}   \right )  \\ 
 \ll \sum_{k \geq 4} e^{k^{1/2}}\exp \left (   -  \kappa  k^{ 1/2}  \log  k    \right )  < \infty\, ,
\end{multline}
where $\kappa$ is a positive constant depending on $a$, $C$ and $p$.  This proves that,  almost surely,
\beq  \label{ine7-13-07}
\lim_{k \rightarrow \infty}  \max_{n_k+1 \leq n \leq n_{k+1}}   \frac{n-n_k}{n}  \sum_{m = 0}^{s_k+1} 2^{pm} \mu(B_m)  \frac{ {\mathcal D}_p ({\mathcal R}_{B_m}\mu_{n_k,n},{\mathcal R}_{B_m} \mu)  }{v_{n_{k+1}}}  =0  \, .
\eeq
Starting from \eqref{ine2-13-07} and taking into account  \eqref{ine5-13-07}, \eqref{ine6-13-07} and  \eqref{ine7-13-07}, it follows that there exists an universal constant $C_p$ depending on $p$ such that 
\[
\limsup_{n \rightarrow \infty}   \frac{ {\mathcal D}_p ( \mu_{n}, \mu)  }{v_{n}}  \leq C_p V   \quad \text{a.s.}
\]
To conclude the case $p >d/2$, it suffices to use inequality  \eqref{lma5FG}.

\medskip

\noindent $\bullet$  If $p\in [1,d/2)$, we proceed as for $p>d/2$, choosing now
\[
v_n = \left ( \frac{\log \log n}{ n }  \right )^{p/d}\, .
\]
Let us give the main steps of the proof. We start again from \eqref{decpreuveLIL}. To deal with the two last terms in the right-hand side of  \eqref{decpreuveLIL}, contrary to the case where $p >d/2$, we do not need here to make a truncation procedure.  Indeed, by Markov's inequality at order $2$, we infer that there exists a positive universal constant $c$ such that, for any $\varepsilon >0$,
\begin{align*}
\P \left (   \max_{n_k \leq n \leq n_{k+1}}   \sum_{m \geq 0} 2^{pm}  | \mu_{n} (B_m) - \mu (B_m) |   
\geq  \varepsilon  V v_{n_{k+1}}  \right )  &  \leq c
   \left ( \frac{ 1 }{  \varepsilon V  n_k v_{n_{k+1}} }  \right )^2  n_{k+1} V^2  \\ 
   & \ll \frac{1}{\varepsilon^{2} n_k^{1-2p/d} (\log \log n_k )^{2p/d}} \, ,
 \end{align*}
which, by an application of the direct part of  Borel-Cantelli's lemma, proves that \eqref{p4-14-dec-18} holds with $\lambda_{\gamma } =0$. Similarly \eqref{p5-14-dec-18} holds and then 
\eqref{ine2-13-07} does also. Hence, it remains to deal with the two last terms in inequality \eqref{ine2-13-07}. This can be done as in the previous case. To handle the probabilities of deviations appearing in \eqref{ine4-13-07bis} and \eqref{ine4-13-07-bisbis}, we proceed again as in the proof of Theorem 2 in \cite{FG} (but this time considering the case $p <d/2$). For instance, concerning the probability of deviation appearing in \eqref{ine4-13-07bis}, this leads to the following inequality: there exists an universal positive constant $a $ such that 
\begin{multline*} 
\P \left (  \sum_{m =0}^{ s_k+1} 2^{pm} \mu(B_m)  {\mathcal D}_p ({\mathcal R}_{B_m}\mu_{n_k},{\mathcal R}_{B_m} \mu)  \geq  C B  v_{n_{k+1}}  \right )  \\
\leq   \sum_{m =0}^{ s_k+1}  \exp \left (   - a C^{d/p} ( \log \log n_k ) \mu(B_m)  \left (  \frac{b_m }{ 2^{pm} \mu ( B_m)} \right )^{d/p}   \right ) 
\, ,
\end{multline*}
where the quantities $s_k$, $B$ and $b_m$ have been introduced previously. 
The probability of deviation appearing in \eqref{ine4-13-07-bisbis} can be handled similarly, and the result follows by taking into account that $(\mu ( B_m) )^{1-p/d} \leq (\mu ( B_m) )^{1/2}$ and inequality \eqref{evident}. 

\subsection{Proof of Theorem \ref{order1}}
Let $q \in (1,2]$ and $M >0$. 
From \eqref{depart} and \eqref{departbis}, we get the upper bound 
$$
 \frac 1 n \Big  \| \max_{1 \leq k \leq n}kW_p^p(\mu_k, \mu)  \Big   \|_1  \leq \frac C n \left (  \Big  \|  \max_{1 \leq k \leq n}kB_{p,M}(\mu_k, \mu)  \Big  \|_1
 +  \Big  \|  \max_{1 \leq k \leq n}kA_{p,M}(\mu_k, \mu)  \Big   \|_q\right ) \, .
$$
Using \eqref{step1}, \eqref{step2} and the same arguments as to get \eqref{Bbound1}, it follows that 
\begin{multline*}
\frac 1 n  \Big    \| \max_{1 \leq k \leq n}kW_p^p(\mu_k, \mu)  \Big   \|_1  \ll \int_0^\infty t^{p-1} H(t) {\bf 1}_{t >M} dt  \\ +
 \sum_{m\geq 0 } 2^{pm} \left( \mu (  B_m  \cap {\mathcal C}_M) \right )^{1/q} \sum_{\ell \geq 0} 2^{-p \ell}   \min \left (1, n^{-(q-1)/q} 2^{\ell d(q-1)/q} \right).
\end{multline*}
Then, using the fact that $ \sum_{m\geq 0 } 2^{pm} ( \mu (  B_m  \cap {\mathcal C}_M) )^{1/q} \ll \int_0^\infty t^{p-1} (H(t))^{1/q} {\bf 1}_{t \leq M} dt$,  we conclude as in Subsection \ref{sec:WVBE} by considering the three cases $p>d(q-1)/q$, $p=d(q-1)/q$ and $p<d(q-1)/q$.

\subsection{Proof of Theorem \ref{order2}}
From  \eqref{depart}, we have that 
$$
 \frac 1 n \left  \| \max_{1 \leq k \leq n}kW_p^p(\mu_k, \mu) \right  \|_2  \leq \frac C n \left  \| \max_{1 \leq k \leq n}k \Delta_p(\mu_k, \mu) \right  \|_2
$$
From \eqref{step1} and \eqref{step2} with $M= \infty$, we get the upper bound 
\begin{align*}
\frac 1 n \left  \| \max_{1 \leq k \leq n}kW_p^p(\mu_k, \mu) \right  \|_2 &  \ll 
 \sum_{m\geq 0 } 2^{pm} \left( \mu (  B_m ) \right )^{1/2} \sum_{\ell \geq 0} 2^{-p \ell}   \min \left (1, 2^{\ell d/2}/\sqrt n \right) \\
 & \ll \left(  \int_0^\infty t^{p-1} (H(t))^{1/2}  dt \right ) \sum_{\ell \geq 0} 2^{-p \ell}   \min \left (1, 2^{\ell d/2}/\sqrt n \right) \, .
\end{align*}
Then we conclude as in Subsection \ref{sec:WVBE} by considering the three cases $p>d/2$, $p=d/2$ and $p<d/2$. 

\subsection{Proof of Theorem \ref{Rosine}}

Let $r >2$. Starting from \eqref{start0}, we infer that, for any positive constant $v_n$, 
\begin{equation}\label{startRos}
\frac{1}{n^r}\left \| \max_{1 \leq k \leq n}kW_p^p(\mu_k, \mu) \right \|_r^r \leq {v_n^r} + r\int_{v_n}^\infty x^{r-1} {\mathbb P}\left ( \max_{1 \leq k \leq n}kW_p^p(\mu_k, \mu) >nx \right ) dx \, ,
\end{equation}
and we  use the upper bound \eqref{start} to deal with the deviation probability in  \eqref{startRos}. Let  $y=x/2C$ and $M >0$. By Markov's inequality at order $q>r$ and 2,
\begin{align}
{\mathbb P}\left ( \max_{1 \leq k \leq n}kA_{p,M}(\mu_k, \mu) >ny\right )  & \leq  \frac{ \left \| \max_{1 \leq k \leq n}kA_{p,M}(\mu_k, \mu)\right \|_q^q}{n^q y^q}  \, , \label{A_pbis}\\
{\mathbb P}\left ( \max_{1 \leq k \leq n}kB_{p,M}(\mu_k, \mu) >ny\right ) & \leq \frac{\left \|  \max_{1 \leq k \leq n}kB_{p,M}(\mu_k, \mu)\right \|_2^2}{n^2 y^2} \label{B_pbis} \, .
\end{align}

To deal with \eqref{B_pbis}, we proceed as to get  \eqref{Abound}, and we obtain
\begin{equation}\label{Bboundbis}
{\mathbb P}\left ( \max_{1 \leq k \leq n}kB_{p,M}(\mu_k, \mu) >ny\right ) \ll 
\frac{1}{y^2} \left ( \sum_{m\geq 0 } 2^{pm} \left( \mu (  B_m \cap  {\mathcal C}^c_M ) \right )^{1/2} \sum_{\ell \geq 0} 2^{-p \ell}   \min \left (1, 2^{\ell d/2}/\sqrt n \right) \right)^2\, .
\end{equation}

Let us now handle  \eqref{A_pbis}. With this aim, in a sake of clarity, let us first recall inequality \eqref{decA-16-dec}, 
\begin{multline*}
\Big \|  \max_{1 \leq k \leq n}kA_{p,M}(\mu_k, \mu)\Big \|_q  \\ \ll 
\sum_{m\geq 0 } 2^{pm} \sum_{\ell \geq 0} 2^{-p \ell} \Big \|  \max_{1 \leq k \leq n}\sum_{ F \in {\mathcal P}_{\ell}} | k\mu_k ( 2^m F \cap B_m \cap  {\mathcal C}_M) - k\mu ( 2^m F \cap B_m \cap  {\mathcal C}_M )  | \Big \|_q \, .
\end{multline*}
By using a maximal version of  Rosenthal's inequality (see for instance \cite{P}),
\begin{multline*}
\Big\| \max_{1 \leq k \leq n}| k\mu_k ( 2^m F \cap B_m \cap  {\mathcal C}_M) - k\mu ( 2^m F \cap B_m \cap  {\mathcal C}_M )  |\Big \|_q
\ll  \sqrt n \left ( \mu ( 2^m F \cap B_m \cap  {\mathcal C}_M ) \right)^{1/2}   \\ 
+ n^{1/q} \left( \mu (  2^m F \cap B_m \cap  {\mathcal C}_M ) \right )^{1/q}\, ,
\end{multline*}
so that, by using H\"older's inequality (twice) and the fact that $| {\mathcal P}_{\ell}|=2^{\ell d}$, 
\begin{multline} \label{todealwithA-16-dec}
\sum_{ F \in {\mathcal P}_{\ell}} \Big \| \max_{1 \leq k \leq n} | k\mu_k (  2^m F \cap B_m \cap  {\mathcal C}_M) - k\mu ( 2^m F \cap B_m \cap  {\mathcal C}_M )  | \Big \|_q 
\leq  2^{\ell d/2} \sqrt{n} \left( \mu (  B_m \cap  {\mathcal C}_M ) \right )^{1/2}  \\ + 2^{\ell d(q-1)/q} n^{1/q} \left( \mu (   B_m \cap  {\mathcal C}_M ) \right )^{1/q}\, .
\end{multline}
So, starting from \eqref{decA-16-dec} and taking into account \eqref{step1} and \eqref{todealwithA-16-dec} together with the fact that for non-negative reals $a, b, c$, 
$\min (a, b+c) \leq \min (a,b) + \min (a,c)$, we get 
\begin{equation}\label{step2bis}
\Big \|  \max_{1 \leq k \leq n}kA_{p,M}(\mu_k, \mu)\Big \|_q  \ll  n(I_1 + I_2) \, ,
\end{equation}
where
$$
I_1=   \sum_{m\geq 0 } 2^{pm} \sum_{\ell \geq 0} 2^{-p \ell} \min \left (\left( \mu (   B_m \cap  {\mathcal C}_M ) \right )^{1/q}, 
n^{-1/2} 2^{\ell d/2}  \left( \mu (  B_m \cap  {\mathcal C}_M ) \right )^{1/2} \right ) 
$$
and 
$$
I_2=  \sum_{m\geq 0 } 2^{pm} \sum_{\ell \geq 0} 2^{-p \ell} \left( \mu (   B_m \cap  {\mathcal C}_M ) \right )^{1/q}  \min \left (1, n^{-(q-1)/q} 2^{\ell d(q-1)/q} \right ) \, .
$$
Combining \eqref{A_pbis} and \eqref{step2bis}, we obtain that 
\begin{equation}\label{Aboundbis}
{\mathbb P}\left ( \max_{1 \leq k \leq n}kA_{p,M}(\mu_k, \mu) >ny\right ) \ll \frac{ (I_1 + I_2)^q}{y^q}  \, .
\end{equation}
From \eqref{Aboundbis}, we see that four cases arise:

\smallskip

\noindent $\bullet$  $p>d(r-1)/r$. In that case $p>d/2$, and 
$$
I_1 \leq  n^{-1/2} \sum_{m\geq 0 } 2^{pm} \left( \mu (  B_m \cap  {\mathcal C}_M ) \right )^{1/2} \ll n^{-1/2} \int_0^\infty t^{p-1} \sqrt{H(t)} {\mathbf 1}_{t  \leq M} dt  \, .
$$
Consequently 
$$
  \int_{v_n}^\infty x^{r-1-q} I_1^q  dx  \ll n^{-q/2} v_n^{r-q}  \left (\int_0^\infty t^{p-1} \sqrt{H(t)}  dt \right )^q  \, .
$$
Choosing $v_n= n^{-1/2} \int_0^\infty t^{p-1} \sqrt{H(t)}  dt $, we get
\begin{equation}\label{firstI1}
\int_{v_n}^\infty x^{r-1-q} I_1^q  dx  \ll n^{-r/2}  \left (\int_0^\infty t^{p-1} \sqrt{H(t)}  dt \right )^r \, .
\end{equation} 

Let us now deal with the term involving $I_2$. First, we choose $q$ close enough to $r$ in such a way that $p>d(q-1)/q$. In that case 
$$
I_2 \leq  n^{-(q-1)/q} \sum_{m\geq 0 } 2^{pm} \left( \mu (  B_m \cap  {\mathcal C}_M ) \right )^{1/q} \ll n^{-(q-1)/q} \int_0^\infty t^{p-1} (H(t))^{1/q} {\mathbf 1}_{t  \leq M} dt  \, .
$$
Let $M=(nx)^{1/p}$. Arguing as in \eqref{Heq} with $\beta < (q-1)/q$, we get 
$$
  \int_0^\infty x^{r-1-q} I_2^q dx \ll \frac{n^{(q-1-\beta q)/p}}{n^{q-1} }\int_0^\infty t^{q(p-1+\beta)} H(t) \int_0^\infty x^{r-1-q}  x^{(q-1-\beta q)/p}{\bf 1}_{x \geq t^p/n} \, dx \, dt  \, .
$$
Taking $\beta$ close enough to $(q-1)/q$ in such a way that $r-q + (q-1-\beta q)/p <0$, we get that
\begin{equation}\label{firstI2}
\int_0^\infty x^{r-1-q} I_2^q dx \ll n^{-(r-1)} \int_0^\infty t^{rp-1} H(t) dt \ll   n^{-(r-1)} \|X\|_{rp}^{rp} Ê\, .
\end{equation} 
From \eqref{Aboundbis}, \eqref{firstI1} and \eqref{firstI2}, we get that 
\begin{equation}\label{part1}
 \int_{v_n}^\infty  x^{r-1} {\mathbb P}\left ( \max_{1 \leq k \leq n}kA_{p,M}(\mu_k, \mu) >nx/(2C)\right )  dx \ll  n^{-r/2}  \left (\int_0^\infty t^{p-1} \sqrt{H(t)}  dt \right )^r + 
 n^{-(r-1)} \|X\|_{rp}^{rp} \, .
\end{equation}

In the same way, since $p>d/2$, we infer from \eqref{Bboundbis} that 
\begin{multline*}
{\mathbb P}\left ( \max_{1 \leq k \leq n}kB_{p,M}(\mu_k, \mu) >ny\right ) \ll \frac{1}{n x^2}  
\left ( \sum_{m\geq 0 } 2^{pm} \left( \mu (  B_m \cap  {\mathcal C}^c_M ) \right )^{1/2} \right )^2 \\
 \ll 
\frac{1}{n x^2}  \left (\int_0^\infty t^{p-1} \sqrt{H(t)} {\mathbf 1}_{t  > M} dt  \right )^2 \, .
\end{multline*}
Proceeding again as in \eqref{Heq} with $\beta > 1/2$, we infer that 
\begin{multline*}
 \int_{v_n}^\infty  x^{r-1} {\mathbb P}\left ( \max_{1 \leq k \leq n}kB_{p,M}(\mu_k, \mu) >nx/(2C)\right )  dx 
 \\
 \ll \frac{n^{(1-2 \beta)/p}}{n} 
 \int_0^\infty t^{2(p-1+\beta)} H(t) \int_0^\infty x^{r-3}  x^{(1- 2 \beta )/p}{\bf 1}_{x < t^p/n} \, dx \, dt  \, .
\end{multline*}
Taking $\beta$ close enough to $1/2$ in such a way that $(r-2)+(1-2 \beta)/p>0$, we get that 
\begin{equation}\label{part2}
 \int_{v_n}^\infty  x^{r-1} {\mathbb P}\left ( \max_{1 \leq k \leq n}kB_{p,M}(\mu_k, \mu) >nx/(2C)\right )  dx 
 \ll n^{-(r-1)} \int_0^\infty t^{rp-1} H(t) dt \ll   n^{-(r-1)} \|X\|_{rp}^{rp} Ê\, .
\end{equation} 

Finally, starting from \eqref{startRos} with  $v_n= n^{-1/2} \int_0^\infty t^{p-1} \sqrt{H(t)}  dt $, and gathering \eqref{start}, \eqref{A_pbis}, \eqref{B_pbis}, \eqref{part1} and \eqref{part2}, Theorem \ref{Rosine} is proved in the case where
$p>d(r-1)/r$.  

\smallskip

\noindent $\bullet$  $ d/2 <p \leq d(r-1)/r$.
In that case we use the upper bound  \eqref{firstI1} without any changes. Let us now deal with the term involving $I_2$.  Starting from the definition of $I_2$, and considering the two cases where either $2^\ell < n^{1/d}$ or $2^\ell \geq n^{1/d}$, we infer that 
$$
I_2 \ll  n^{-p/d} \sum_{m\geq 0 } 2^{pm} \left( \mu (  B_m \cap  {\mathcal C}_M ) \right )^{1/q} \ll n^{-p/d} \int_0^\infty t^{p-1} (H(t))^{1/q} {\mathbf 1}_{t  \leq M} dt  \, .
$$
Let $M=(nx)^{1/p}/u_n$ for some sequence of positive numbers $(u_n)_{n >0}$ that will be chosen later. Arguing as in \eqref{firstI2}, we get 
\begin{equation}\label{secondI2}
\int_0^\infty x^{r-1-q} I_2^q dx \ll n^{q-r-pq/d} u_n^{p(r-q)}  \int_0^\infty t^{rp-1} H(t) dt \ll   n^{q-r-pq/d} u_n^{p(r-q)} \|X\|_{rp}^{rp} Ê\, .
\end{equation} 
In the same way, arguing as to get \eqref{part2},
\begin{multline}\label{part2bis}
 \int_{v_n}^\infty  x^{r-1} {\mathbb P}\left ( \max_{1 \leq k \leq n}kB_{p,M}(\mu_k, \mu) >nx/(2C)\right )  dx  \\
 \ll n^{-(r-1)}  u_n^{p(r-2)} \int_0^\infty t^{rp-1} H(t) dt \ll   n^{-(r-1)} u_n^{p(r-2)} \|X\|_{rp}^{rp} Ê\, .
\end{multline} 
Now $n^{ q-r-pq/d }  u_n ^{(r - q)p } = n^{-(r-1)} u_n^{p(r-2)}$  iff  $u_n ^{p } = n^{-1/(q-2)} n^{ (1 -p /d) q / (q-2) } $. With this choice of $u_n$ and taking $q = r + \varepsilon$, we have 
\[
n^{ q-r-pq/d }  u_n ^{(r - q)p } = n^{-r p/d} (n^{(d-p)/d}/u_n^p)^{q-r} = n^{-r p/d} n^{(2p-d)(q-r)/(d(q-2))}= n^{-r p/d} n^{\varepsilon (2p-d)/(d(r-2 + \varepsilon))} \, .
\]
Hence, with this choice of $u_n$, the upper bounds  \eqref{firstI1}, \eqref{secondI2} and \eqref{part2bis} give the desired inequality for 
$ d/2 <p \leq d(r-1)/r$. 

\smallskip

\noindent $\bullet$  $ p < d/2$.
Note first that, by homogeneity, the general inequality may be deduced from the case where $\|X\|_{rp}=1$ by considering the variables 
$X_i/\|X\|_{rp}$. Hence, from now, we shall assume that $\|X\|_{rp}=1$. 

Let $M=(nx)^{1/p}/u_n$ for some sequence of positive numbers $(u_n)_{n >0}$. We first note that, since $q>d/(d-p)$ (indeed $q >2$ and $d/(d-p) < 2$), the upper bound \eqref{secondI2} holds.
Taking $u_n= n^{1/p}/n^{1/d}$, we get 
\begin{equation}\label{thirdI2}
\int_0^\infty x^{r-1-q} I_2^q dx  \ll   n^{-rp/d} \, .
\end{equation} 

 Let us now deal with the term involving $I_1$.
Starting from the definition of $I_1$, and considering the two cases where 
\[ \text{ either } \quad 2^\ell < n^{1/d} (\mu (  B_m \cap  {\mathcal C}_M ))^{(2-q)/(dq)} \quad  \text{or}  \quad 2^\ell \geq n^{1/d}  (\mu (  B_m \cap  {\mathcal C}_M ))^{(2-q)/(dq)} \, , \]
we infer that 
\begin{equation}\label{interm}
I_1 \ll  n^{-p/d} \sum_{m\geq 0 } 2^{pm} \left( \mu (  B_m \cap  {\mathcal C}_M ) \right )^{(d+p(q-2))/(dq)} \ll n^{-p/d} \int_0^\infty t^{p-1} (H(t))^{(d+p(q-2))/(dq)}  dt  \, .
\end{equation}
We choose now $q>r$ such that $(d+p(q-2))/(dq)>1/r$ (this is true whatever $q$ if $p\geq d/r$, otherwise we need to choose 
$r<q< r(d-2p)/(d-rp)$). Since $\|X\|_{rp}=1$,  $H(t) \leq \min(1, t^{-rp})$, which together \eqref{interm}  and the choice of $q$ implies that 
$I_1 \ll n^{-p/d}$. Consequently, taking $v_n=n^{-p/d}$,
\begin{equation}\label{thirdI1}
  \int_{v_n}^\infty x^{r-1-q} I_1^q \, dx \ll n^{-qp/d} v_n^{r-q}  \ll n^{-rp/d} \, .
\end{equation}
From \eqref{Aboundbis}, \eqref{thirdI2} and \eqref{thirdI1}, we get that 
\begin{equation}\label{part1third}
 \int_{v_n}^\infty  x^{r-1} {\mathbb P}\left ( \max_{1 \leq k \leq n}kA_{p,M}(\mu_k, \mu) >nx/(2C)\right )  dx \ll  n^{-rp/d}  \, .
\end{equation}

On another hand, since $p<d/2$, we infer from \eqref{Bboundbis} that 
\begin{multline*}
{\mathbb P}\left ( \max_{1 \leq k \leq n}kB_{p,M}(\mu_k, \mu) >ny\right ) \ll \frac{1}{n^{2p/d} x^2}  
\left ( \sum_{m\geq 0 } 2^{pm} \left( \mu (  B_m \cap  {\mathcal C}^c_M ) \right )^{1/2} \right )^2 \\
 \ll 
\frac{1}{n^{2p/d} x^2}  \left (\int_0^\infty t^{p-1} \sqrt{H(t)} {\mathbf 1}_{t  > M} dt  \right )^2 \, .
\end{multline*}
Proceeding again as in \eqref{part2bis}, we get
\begin{equation}\label{part2third}
 \int_{v_n}^\infty  x^{r-1} {\mathbb P}\left ( \max_{1 \leq k \leq n}kB_{p,M}(\mu_k, \mu) >nx/(2C)\right )  dx  
 \ll \frac{u_n^{p(r-2)}}{ n^{((r-2)d+2p)/d }} \int_0^\infty t^{rp-1} H(t) dt \ll  n^{-rp/d}  Ê\, ,
\end{equation} 
the last inequality being true because $u_n= n^{1/p}/n^{1/d}$ and $\|X\|_{rp}=1$. 

Finally, starting from \eqref{startRos} with  $v_n =n^{-p/d}$, and gathering \eqref{start}, \eqref{part1third} and \eqref{part2third}, Theorem \ref{Rosine} is proved in the case where
$p< d/2 $ and $\|X\|_{rp}=1$.

\smallskip

\noindent $\bullet$  $ p = d/2$. Again, without loss of generality we can assume that $\|X\|_{rd/2}=1$. 
We proceed as before to handle the term $\int_{v_n}^{\infty} x^{r-1} {\mathbb P} ( \max_{1 \leq k \leq n}kA_{p,M}(\mu_k, \mu)  > nx/(2C)  ) dx$. We take $q >r$ and use the Rosenthal inequality. We then infer that 
\[
I_1 \ll n^{-1/2}  \log n  \left (  \int_0^{\infty} t^{d/2 -1} \sqrt{H(t)} {\bf 1}_{t \leq M } dt  \right )  +  n^{-1/2}   \left (  \int_0^{\infty} t^{d/2 -1} \sqrt{H(t)} \log (1/H(t) ) {\bf 1}_{t \leq M } dt  \right )  \, .
\]
Therefore, if we choose
\begin{equation} \label{contrainte1vn}
v_n \geq n^{-1/2} \max \left (  \log n  \int_0^{\infty} t^{d/2 -1} \sqrt{H(t)} dt ,  \int_0^{\infty} t^{d/2 -1} \sqrt{H(t)} \log (1/H(t) )  dt  \right ) =: v_n(1) \, ,
\end{equation}
we get 
\[
 \int_{v_n}^{\infty} x^{r-1-q} I_1^q dx \ll n^{-r/2} ( \log n )^r  \left (  \int_0^{\infty} t^{d/2 -1} \sqrt{H(t)} dt  \right )^r + n^{- r/2}  \left (  \int_0^{\infty} t^{d/2 -1} \sqrt{H(t)} \log (1/H(t) )  dt 
  \right )^r \, .
\]
Since $H(t) \leq \min (1, t^{-rd/2} ) $ and $r>2$, it follows that 
\begin{equation}\label{I1d/2}
 \int_{v_n}^{\infty} x^{r-1-q} I_1^q dx \ll n^{-r/2} ( \log n )^r  \left (  \int_0^{\infty} t^{d/2 -1} \sqrt{H(t)} dt  \right )^r + n^{- r/2}   \, .
\end{equation}
On another hand, we have
\[
 \int_{v_n}^{\infty} x^{r-1-q} I_2^q dx \ll n^{-q/2} \int_{v_n}^{\infty} x^{r-1-q} \left (  \int_0^{\infty} t^{d/2 -1} H^{1/q}(t)  {\bf 1}_{t \leq M } dt  \right )^q dx  \, .
\]
Selecting
\[
M = (nx)^{2/d} /u_n \text{ with } u_n =n^{1/d} \, , 
\]
we get, by taking into account previous computations, that 
\begin{equation}\label{I2d/2}
\int_{v_n}^{\infty} x^{r-1-q} I_2^q dx \ll n^{-r/2} \int_0^{\infty} t^{ rd/2 -1 } H (t) dt = n^{-r/2}  \, .
\end{equation}

We handle now the quantity  $\int_{v_n}^{\infty} x^{r-1} {\mathbb P} ( \max_{1 \leq k \leq n}kB_{p,M}(\mu_k, \mu)  > nx/(2C)  ) dx$.  We shall apply this time the Rosenthal inequality as we did to handle $ \Vert   \max_{1 \leq k \leq n}kA_{p,M}(\mu_k, \mu)   \Vert_q $, but with $q \in (2,r)$.  We then infer that 
\begin{equation}\label{Bpd/2}
\left  \Vert   \max_{1 \leq k \leq n}kB_{p,M}(\mu_k, \mu)  \right  \Vert_q \ll  n ( J_1 + J_2  + J_3) \, ,
\end{equation}
with 
\[
J_1 = n^{-1/2}  \log n    \int_0^{\infty} t^{d/2 -1} \sqrt{H(t)} {\bf 1}_{t >M } dt   \, , 
\]
\[
J_2 =   n^{-1/2}   \int_0^{\infty} t^{d/2 -1} \sqrt{H(t)} \log (1/H(t) ) {\bf 1}_{t> M } dt   \, , 
\]
and 
\[
J_3 =   n^{-1/2}   \int_0^{\infty} t^{d/2 -1} H^{1/q}(t)  {\bf 1}_{t> M } dt   \, . 
\]
Note that since $M = (nx)^{2/d} /u_n $ with $u_n =n^{1/d}$, applying   H\"older's inequality as in previous computations, we get
\begin{equation}\label{J3d/2}
 \int_{v_n}^{\infty} x^{r-1-q} J_3^q dx  \ll n^{-r/2} \int_0^{\infty} t^{ rd/2 -1 } H (t) dt  \, .
\end{equation}
On another hand, using that $H(t) \leq \min (1, t^{-rd/2} ) $, we have (since $r>2$ and  $M^{d/2} = x \sqrt{n}$), 
\begin{multline*}
 \int_{v_n}^{\infty} x^{r-1-q} J_1^q dx  \leq  n^{-q/2} (\log n)^q  \int_{v_n}^{\infty} x^{r-1-q}  \left (  \int_M^{\infty} t^{d/2 -1} t^{-rd/4} dt \right )^q dx 
 \\ \ll  n^{-rq/4} (\log n)^q  \int_{v_n}^{\infty} x^{r (1-q/2)-1}  dx  \ll n^{-rq/4} (\log n)^q  v_n^r v_n^{-rq/2}  \, .
\end{multline*}
Therefore if 
\begin{equation} \label{contrainte2vn}
v_n \geq  n^{-1/2}  (  \log n )^{2/r} =:v_n(2) \, ,
\end{equation}
we get 
\begin{equation} \label{J1d/2}
\int_{v_n}^{\infty} x^{r-1-q} J_1^q dx   \ll n^{-r/2}   (  \log n )^{2}   \, .
\end{equation}
We handle now the term involving $J_2$. We have
\[  
 \int_{v_n}^{\infty} x^{r-1-q} J_2^q dx  =   n^{-q/2} \int_{v_n}^{\infty} x^{r-1-q} \left (  \int_0^{\infty} t^{d/2 -1} \sqrt{H(t)} \log (1/H(t) ) {\bf 1}_{t> M } dt \right )^q dx  \, .
\]
If $v_n \geq n^{-1/2}$, using that $H(t) \leq \min (1, t^{-rd/2} ) $, simple computations lead to 
\begin{equation*}
  \int_{v_n}^{\infty} x^{r-1-q} J_2^q dx   
\ll   n^{-r/2}    (\sqrt{ n}v_n)^{r(1 -q/2)}   \{  (\log ( \sqrt {n} v_n ) )^{q}  +1 \}  \, .
\end{equation*}
Therefore, if \eqref{contrainte2vn} holds, 
we get
\begin{equation} \label{J2d/2}
 \int_{v_n}^{\infty} x^{r-1-q} J_2^q dx  \ll   n^{-r/2}      (\log n )^2    \, .
\end{equation}
So finally if we choose
\[
v_n = \max (v_n (1), v_n(2) ) \, ,
\]
the constraints \eqref{contrainte1vn} and  \eqref{contrainte2vn}  are satisfied. 
Starting from \eqref{startRos}, and gathering  the bounds \eqref{start}, \eqref{I1d/2}, \eqref{I2d/2}, \eqref{Bpd/2},  \eqref{J3d/2}, \eqref{J1d/2},  and \eqref{J2d/2},
 we get the desired inequality in the  case  $\|X\|_{rd/2}=1$. 

\bigskip

\bigskip

\noindent {\bf  \large Acknowledgements:} We thank the referee for his/her precise comments, and for pointing out some important references that were missing in the previous version
of this paper.

\end{document}